\documentclass[11pt,letterpaper]{article} 

\usepackage{amssymb,latexsym,amsmath} 
\usepackage[dvipsnames]{color} 
\usepackage{graphics} 

\definecolor{cite}{named}{Mahogany} 
\definecolor{link}{named}{OliveGreen} 
\definecolor{anchor}{named}{Gray} 

\def\QED{ 
\hfill%\scalebox{.04}{\includegraphics{fondoe}}} 
$\Box$} 
\def\demo{\noindent{\bf Proof. }} 

\def\NN{\mathbb{N}}

\def\b1{{\bf 1}}

\def\R_A{{\mathbb R}_+{\cal A}} 
 
\def\min{\rm min}

\def\deg{\mathop{\rm deg}\nolimits}

\newcommand{\bin}[2]{ 
   \left( 
     \begin{array}{@{}c@{}} 
         #1  \\ #2 
     \end{array} 
   \right)          }

\newtheorem{Theorem}{Theorem}[section] 
\newtheorem{Definition}[Theorem]{Definition} 
\newtheorem{Lemma}[Theorem]{Lemma} 
\newtheorem{Corollary}[Theorem]{Corollary} 
\newtheorem{Conjecture}[Theorem]{Conjecture} 
\newtheorem{Remark}[Theorem]{Remark} 
\newtheorem{Example}[Theorem]{Example\/} 
 
\newtheorem{Proposition}[Theorem]{Proposition} 
\newtheorem{Claim}[Theorem]{Claim}

\newcounter{aux1}
\newcounter{aux2}

\begin{document} 
\topmargin3mm 

\medskip 

%===================================% 
%==============Titulo===============% 
%===================================% 

\vspace{1.5cm} 
\begin{center} 
{\Large\bf On bounds for some graph invariants} 
\\ 
\vspace{3mm} 
Isidoro Gitler and Carlos E. Valencia
\footnote{The authors where partially supported by CONACyT grant 49835 and SNI.}     
\\ {\small Departamento de Matem\'aticas} 
\\ {\small Centro de Investigaci\'on y de Estudios Avanzados del IPN} 
\\ {\small Apartado Postal 14--740} 
\\ {\small 07000 M\'exico City, D.F.} 
\\ {\small e-mail: {\tt igitler@math.cinvestav.mx}, {\tt cvalencia75@gmail.com}} 
\\ 
\medskip 

\end{center}

%===================================% 
%==============Resumen==============% 
%===================================% 

\begin{abstract} 
\noindent 
Let $G$ be a graph without isolated vertices, $\alpha(G)$ and $\tau(G)$
be the stability number and the covering number of $G$, respectively. 

The paper is divided in two parts:
In the first part we study the minimum number of edges that a $k$-connected
graph can have as a function of $\alpha(G)$ and $\tau(G)$. 
In particular, we obtain the following lower bound: 
$$ 
q(G) \geq \alpha(G)-c(G)+ \Gamma(\alpha(G), \tau(G)), 
$$ 
where $q(G)$ is the number of edges of $G$, $c(G)$ is the number of connected components of $G$ and
$$
\Gamma(\alpha(G),\tau(G))=(\alpha(G)-s)  \bin{r}{2}+s  \bin{r+1}{2},
$$ 
where $\alpha(G)+\tau(G)=r\alpha(G)+s$ with $0\leq s< \alpha(G)$. 

This is a solution to an open question posed by Ore in his 
book~\cite[pag. 216]{Ore62}, which indeed is a variant for connected graphs 
of a celebrated theorem of Tur\'an~\cite{turan41}. 

In the second part of this paper, we study the relations between 
$\alpha(G)$, $\tau(G)$ and $\delta(G) = \alpha(G)-\sigma_{v}(G)$, 
where the {\it $\sigma_{v}$-cover} number of a graph, denoted by $\sigma_{v}(G)$, 
is the maximum natural number $m$, such that every vertex of 
$G$ belongs to a maximal independent set with at least $m$ vertices.   
The main theorem of this part states that 
\[
\alpha(G)\leq \tau(G)[1+\delta(G)].
\]  

In the last section, we discuss some conjectures related to this theorem. 

\end{abstract} 

%==================================================================% 
%===========================Introducion============================% 
%==================================================================% 

\section{Introduction}

Given a graph $G=(V(G),E(G))$, a subset $M\subseteq V(G)$ is a {\it stable} set if no two vertices in $M$ are adjacent.  
We say that $M$ is a {\it maximal stable} set if it is maximal with respect to inclusion. 
The {\it stability number} of $G$ is given by 
$$ 
\alpha(G)={\rm max}\{|M| \, | \, M\subset V(G) \mbox{ is a stable set in } G\}. 
$$
Also, $C\subseteq V(G)$ is a {\it vertex cover} for a graph
$G$ if every edge of $G$ is incident with at least one vertex in $C$. 
Moreover, the vertex cover $C$ is called a {\it minimal vertex cover} 
if there is no proper subset of $C$ which is a vertex cover.
It is convenient to regard the empty set as a minimal vertex cover 
for a graph with all its vertices isolated. 

The {\it vertex covering number} of $G$, denoted by $ \tau (G)$, 
is the number of vertices in a minimum vertex cover in $G$, that is, 
the size of any smallest vertex cover in $G$. 
Note that a set of vertices in $G$ is a maximal stable set if and only if its 
complement is a minimal vertex cover for $G$, thus $\alpha(G)+\tau(G)=|V(G)|$.

%In this paper, we prove some inequalities relating the stability number 
%$\alpha(G)$ %(maximum size of an independent or stable set of the graph G) 
%and the covering number $\tau(G)$ %(minimum size of a vertex cover of G) 
%of a graph $G$. 

This paper has two main parts. 
In the first part we solve a problem posed by Ore in his book~[pag. 216, research prob. 1]\cite{Ore62}, whose statement is:
\medskip

{\it 
Determine the connected graphs satisfying  $\alpha(G) <k$ ($3\leq k\leq |V(G)|$) and having a minimal number of edges.
} 

\medskip

%What is the minimum number of edges that a connected graph on n vertices and stability number can have? 
This question was completely solved by Tur\'an~\cite{turan41} (see also \cite[pages 214--216]{Ore62}) when the word connected is removed. 
%Its solution is the dual of the famousTuran's Theorem~\cite{turan41}. 
%In this case the optimal graphs are the disjoint union of complete graphs, 
%more precisely the complements of Tur\'an's graphs. 
In the connected case, we completely solve the question by giving a tight lower bound on the number 
of edges of a graph of a given order and with given stability number
(Theorem~\ref{minimal} for connected graphs). 
This result was independently proved by J. Christophe et al in \cite{doignon}.
We include the full classification of all optimal graphs (those achieving the bound). 
A preliminary version of theorems~\ref{minimal} and~\ref{alon} and lemma~\ref{conexo} appear by first time in~\cite{BoundsIIArxiv} and ~\cite{BoundsIIMor}.
The corresponding result for $2$-connected graphs together
with a full classification of the optimal graphs is also included (Theorem~\ref{2conexo}).
%After that this paper was submitted, the autho.
The classification of the connected and $2$-connected optimal graphs was obtained independently in~\cite{bougard}.

\medskip

In the second part of the paper, we prove the inequality $\alpha(G) \leq \tau(G) (1 + \delta(G))$ where
$\delta(G) = \alpha(G)- \sigma_v(G)$ and $\sigma_v(G)$ is the largest $m$, such that every vertex of G belongs to an
independent set of size $m$. As a corollary we obtain $\alpha(G)-|\alpha_{core}| \leq  \tau(G)- | \tau_{core}|$, 
where the  $\alpha_{core}$ and the $\tau_{core}$ are defined as the intersection of all maximum-size
stable sets and the intersection of all minimum-size vertex coverings, respectively. 

The origin of our interest in the study of these relations comes from monomial algebras. 
More precisely, the stability number $\alpha(G)$ of a graph $G$, is equal to
the dimension of the Stanley-Reisner ring associated to the graph $G$;
and the covering number $\tau(G)$ of $G$ is equal to the height
of the ideal associated to the graph $G$. 
Finally, $\sigma_v(G)$ is an upper bound to the depth of this ring. 

From the algebraic point of view, an important class of rings 
is given by those rings $R$ such that their dimension is equal to their depth.
The rings in this class are called Cohen-Macaulay rings.
A graph is Cohen-Macaulay if the Stanley-Reisner ring associated to it, is Cohen-Macaulay.
If a graph $G$ is Cohen-Macaulay, then $\delta(G)=0$ (\cite[Proposition 6.1.21]{monalg}).
Note that this is a necessary, but not a sufficient condition.
The family of graphs with $\delta(G) \geq 1$ corresponds to the 
Stanley-Reisner rings that have a large depth. 
Moreover, the dimension minus the depth is bounded below by $\delta(G)$, and hence $\delta(G)$ is a measure of 
how far these rings are from being Cohen-Macaulay.

\medskip

The outline of the article is as follows: 
We begin with 
%Section 2 of preliminaries, where we briefly introduce the main subject 
%background and basic concepts related to the problems studied in the paper. 
%
%We properly begin the exposition of our results in 
Section 2, where we solve the low 
connectivity (one and two connected) versions of Tur\'an's theorem as thoroughly explained above. 
In this section we give a lower bound for the number of edges of a graph
as a function of its stability and covering numbers (Theorem~\ref{minimal}) together with a characterization of 
$q$-minimal (Lemma~\ref{conexo}) %(connected) 
and $\{q,2\}$-minimal ($2$-connected) graphs (Theorem~\ref{2conexo}). 
%we determine the value of $f (n, a, 2)$ and characterize (n, a, 2)-extremal graphs. 

In section~\ref{bounds} we study some relations between the 
stability and covering numbers of a graph.
Specifically we prove the main result (Theorem~\ref{alon}) of this section,
which is an inequality that gives an upper bound for the stability number of a graph 
with respect to the covering number and $\delta(G) = \alpha(G)- \sigma_v(G)$.
This result generalizes an inequality given in \cite{bounds} which was only valid for $B$-graphs. 
Then we introduce the $\alpha_{core}$ and the $\tau_{core}$ of a graph and relate them by an 
inequality with the stability and covering numbers of the graph. 
Finally, we give a series of conjectures that relate several invariants of graphs for 
$B$-graph and hypergraphs. 

%We change of focus and in section~\ref{bounds} we study some relations between the stability 
%and covering numbers of a graph to obtain some interesting bounds as a consequence of these results.
%Specifically we prove Theorem~\ref{alon}, which is the main theorem of this section. 
%This result is an inequality giving an upper bound for the stability number of a graph 
%with respect to the covering number and $\delta(G) = \alpha(G)- \sigma_v(G)$, where 
%$\sigma_v(G)$ is the largest m such that every vertex of $G$ belongs to an independent set of size $m$. 
%This result generalizes an inequality given in \ref{bounds} which was only valid for $B$-graphs. 
%Then we introduce the $\alpha_{core}$ and the $\tau_{core}$ of a graph and relate them by an 
%inequality with the stability and covering number of the graph. Finally, we give a series of 
%conjectures that relate several invariants of graphs for $B$-graph and hypergraphs. 

%==================================================================% 
%==========================Preliminaries===========================% 
%==================================================================% 

%\section{Preliminaries} 

In this article, all graphs are supposed to be finite and simple (i.e., without loops
and multiple edges).
Let $G=(V, E)$ be a graph with $|V|=n$ vertices and $|E|=q$ edges. 
Given a subset $U\subset V$, the {\it neighbour set} of $U$, 
denoted by $N(U)$, is defined as 
$ N(U)=\{v\in V\mid v \; \mbox{\rm is adjacent to some vertex in}\; U\}$. 

A subset $W$ of $V$ is called a {\it clique} if any two vertices in $W$ are adjacent. 
We call $W$ {\it maximal} if it is maximal with respect to inclusion. 
The {\it clique number} of a graph $G$ is given by 
$$ 
\omega(G)={\rm max}\{ |W| \, | \, W\subset V(G) \mbox{ is a clique in }  G\}. 
$$ 

The {\it complement} of a graph $G$, denoted by $\overline{G}$, is the graph 
with the same vertex set as $G$, and edges all pairs of distinct vertices that 
are nonadjacent in $G$. 
Clearly, $W$ is a clique of $G$ if and only if $W$ is a stable set 
of $\overline{G}$, and thus $\omega(G)=\alpha(\overline{G})$. 

A subgraph $H$ is called an {\it induced subgraph} of $G$, 
denoted by $G[V(H)]$, if $H$ contains all the edges 
$\{v_i,v_j\}\in E(G)$ with $v_i,v_j\in V(H)$.

%A graph $G$ is called {\it $k$-connected} if it has $k$ vertex disjoint paths between 
%every pair of distinct vertices $u$ and $v$.
 
A non-empty graph $G$ is called connected if any two of its vertices are linked by a path in $G$. 
A graph $G$ is called {\it $k$-connected}  (for $k\in {\mathbb N}$) if $|V(G)|>k$ and
$G\setminus X$ is connected for every set $X\subseteq V(G)$ with $|X|<k$.

%==================================================================% 
%============================Seccion 1=============================% 
%==================================================================% 

\section{Low connectivity versions of Tur\'an's theorem}
%Graphs with a minimal number of edges} 

In this section, we study the minimal number of edges in $k$-connected graphs.
Theorem~\ref{minimal} establishes a lower bound for the number of edges of a graph $G$ 
as a function of $\alpha(G)$, $\tau(G)$ and its number of connected components, $c(G)$. 
%This is an answer to an open question posed by Ore in his book~\cite{Ore62}, which is 
%a variant for connected graphs of a celebrated theorem of Tur\'an~\cite{turan41}.  
As a byproduct of the proof of Theorem~\ref{minimal}, we find a 
bound for $2$-connected graphs and determine the graphs for which these bounds are sharp.

A {\it Tur\'an graph}, denoted by $T(a, t)$, is a graph of order $a+t$
consisting of the disjoint union of $a-s$ cliques of order $r= \lfloor \frac{a+t}{a} \rfloor$ and $s$ cliques 
of order $r+1$, where $a+t= ra+s$ with $0 \leq s < a$.

For a graph $G=(V,E)$, we  denote by $q(G)$ the cardinality of its edge set $E(G)$. 
We say that a $k$-connected graph $G$ is {\it $\{q,k\}$-minimal}, if 
there is no graph $G'$ such that 
\begin{description} 
\item[{\it (i)}] $G'$ is $k$-connected, 
\item[{\it (ii)}] $\alpha(G')=\alpha(G)$, 
\item[{\it (iii)}] $\tau(G')=\tau(G)$, and 
\item[{\it (iv)}] $q(G')<q(G)$. 
\end{description} 

We say that an edge $e$ of a graph $G$ is an {\it $\alpha$-critical edge} if $\alpha(G-e)=\alpha(G)+1$.
A vertex $v$ of a graph $G$ is a {\it $\tau$-critical vertex} if $\tau(G-v)=\tau(G)-1$.
A connected graph $G$ is called {\it $\alpha$-critical} ({\it $\tau$-critical}) if all its edges 
(vertices) are $\alpha$-critical ($\tau$-critical).
In Chapter 12 of the book of Lovasz and Plummer~\cite{LovaszPlummer} some of the basic properties of $\alpha$-critical graphs can be found.
For instance, Corollary 12.1.8 in \cite{LovaszPlummer} says that every $\alpha$-critical graph is $2$-connected.
Also, by Lemma 12.1.2 in \cite{LovaszPlummer}, if $G$ is an $\alpha$-critical graph without isolated vertices, 
then $\alpha(G)=\alpha(G-v)$ for all $v\in V(G)$.
Using the previous observation and the fact that $\alpha(G)+\tau(G)=|V(G)|$ we can conclude that if 
$G$ is an $\alpha$-critical graph, then $G$ is a $\tau$-critical graph.

\medskip

For simplicity a $\{q,1\}$-minimal graph will be called a $q$-minimal graph.
Hence, if $G$ is $q$-minimal, then either $\alpha(G)<\alpha(G-e)$ or $c(G)<c(G-e)$
for all the edges $e$ of $G$ (note that $\alpha(G)<\alpha(G-e)$ if and only if $\tau(G)>\tau(G-e)$).   
That is, an edge of a $q$-minimal graph is either $\alpha$-critical or a bridge. 
Therefore the blocks (a maximal connected subgraph without a cutvertex) of a $q$-minimal graph are $\alpha$-critical graphs. %A {\it block} of a graph is 
%Here a graph is $\alpha$-critical if $\alpha(G-e)=\alpha(G)+1$ 
%for all the edges $e$ of $G$ and it is a $\tau$-critical graph if 
%$\tau(G-v)=\tau(G)-1$ for all the vertices $v$ of $G$.

\medskip 

In order to bound the number of edges of a graph we introduce the following numerical function. 
For any natural numbers $a$ and $t$, let 
\[
\Gamma(a,t)=(a-s)  \bin{r}{2}+s  \bin{r+1}{2},
\]
where $a+t=ra+s$ with $0\leq s< a$. 
In other words, $r=1+\left \lfloor \frac{t}{a} \right \rfloor$ and $s=t+a \left \lfloor \frac{t}{a}\right \rfloor$.

%%%%%%%%%%%%%%%%%%%%%%%%%%%%%%%%%%%%%%%%%%%%%%%%%%%%%%%%%%%%%%%%%%%%%%%%%%%%%%%%%%%%%%%
%%%%%%%%%%%%%%%%%%%%%%%%%%%%Lemma Propiedades de Gamma%%%%%%%%%%%%%%%%%%%%%%%%%%%%%%%%%
%%%%%%%%%%%%%%%%%%%%%%%%%%%%%%%%%%%%%%%%%%%%%%%%%%%%%%%%%%%%%%%%%%%%%%%%%%%%%%%%%%%%%%%

\begin{Lemma}\label{lemaG} \setcounter{aux1}{\value{section}} \setcounter{aux2}{\value{Theorem}}
Let $a$ and $t$ be natural numbers, then 
\begin{description} 
\item[{\it (i)}] $ \Gamma(a,t)= 
{\rm min} \left\{ \sum_{i=1}^{a}\bin{z_i}{2}\, : \, z_1 + \cdots + z_{a}= 
a+t \mbox{ and } z_i \geq 0 \mbox{ for all }  1 \leq i \leq a \right\}$ .

\item[{\it (ii)}] $\Gamma(a-1,t)-\Gamma(a,t) \geq \frac{1}{2}(\left \lfloor \frac{t}{a}\right \rfloor^2-\left \lfloor \frac{t}{a}\right \rfloor)=
\bin{\left \lfloor \frac{t}{a}\right \rfloor+1}{2} \geq 0$ 
for all $a\geq 2$ and $t \geq 1$. 

Moreover, $\Gamma(a-1,t)-\Gamma(a,t) = \bin{\left \lfloor \frac{t}{a}\right \rfloor+1}{2}$
if and only if $1+\left \lfloor \frac{t}{a}\right \rfloor \geq \frac{t}{a-1}$ 
and $\bin{\left \lfloor \frac{t}{a}\right \rfloor+1}{2}=0$ if and only if $0\leq t < a$. 

\item[{\it (iii)}] $\Gamma(a,t)-\Gamma(a,t-1)=1+\left \lfloor \frac{t-1}{a}\right \rfloor=\left \lceil \frac{t}{a}\right \rceil$ 
for all $a\geq 1$ and $t \geq 2$.

\item[{\it (iv)}] $\sum_{i=1}^k \Gamma(a_i,t_i) \geq \Gamma(\sum_{i=1}^k a_i, \sum_{i=1}^k t_i)$ 
for all $a_i\geq 1$ and $t_i \geq 1$. 

Furthermore, if $a_1, a_2 \geq 2$, then
$$ 
\Gamma(a_1,t_1)+\Gamma(a_2,t_2)=\Gamma(a_1+a_2,t_1+t_2) 
$$ 
if and only if either $\left \lfloor \frac{t_1}{a_1}\right \rfloor=\left \lfloor \frac{t_2}{a_2}\right \rfloor$,
$\left \lfloor \frac{t_1}{a_1}\right \rfloor-\left \lfloor \frac{t_2}{a_2}\right \rfloor=1$ and $t_1=r_1a_1$ 
or $\left \lfloor \frac{t_2}{a_2}\right \rfloor-\left \lfloor \frac{t_1}{a_1}\right \rfloor=1$ and $t_2=r_2a_2$. 

\item[{\it (v)}] $\left \lceil \frac{2(a-1+\Gamma(a,t))}{a+t} \right \rceil 
= 1 +\left \lfloor \frac{t}{a}\right \rfloor+L$, where $-1 \leq L \leq 1$.
 
Moreover, $L=-1$ if and only if $a=1$. %, $L=1$ if $1 \leq t <a$ or $a+2< t< 2a$ and $L=0$ otherwise.
\end{description} 
\end{Lemma}
\demo 
{\bf (i)} For $a=1$ the result is trivial. For $a \geq 2$ we  use the next observation: 
Let $n,m \geq 1$ be natural numbers with $n > m+1$, then 
\[ 
\bin{n}{2}+\bin{m}{2} > \bin{n-1}{2}+\bin{m+1}{2}. 
\] 

Let $a \geq 2$ and $t \geq 1$ be fixed natural numbers, 
$(z_1,\ldots,z_{a})\in \NN^{a}$ such that $\sum_{i=1}^{a} z_i=a+t$ and let 
$L(z_1,\ldots,z_{a})=\sum_{i=1}^{a}\bin{z_i}{2}$. 
Now, if 
$$ 
\{z_1,\ldots,z_{a}\}\neq \{ \underbrace{r,\ldots,r}_{a-s},\underbrace{r+1,\ldots,r+1}_{s}\} 
$$ 
where $a+t=ra+s$ with $0\leq s< a$, 
then there exist $z_{i_1}$ and $z_{i_2}$ with $ z_{i_1} > z_{i_2}+1$. 
Applying the previous observation we obtain that 
$$ 
L(z_1,\ldots z_{a}) > L(z_1,\ldots,z_{i_1}-1,\ldots,z_{i_2}+1,\ldots,z_{a}) \geq \Gamma(a,t), 
$$ 
and therefore we obtain the result. 

{\bf (ii)} 
Let $a+t=a r+s$ with %$r \geq 1$ and 
$0\leq s< a$, 
then 
$$ 
a+t-1=(a-1)(r+l)+s' 
$$ 
where $r+s-1=(a-1)l+s'$ with $l\geq 0$ and $0\leq s'< a-1$. 

%Using part (i) and 
After some algebraic manipulations we obtain that 
$$ 
2(\Gamma(a-1,t)-\Gamma(a,t))=(r^2-r)+(l^2-l)(a-1)+2ls'. 
$$   
Therefore $\Gamma(a-1,t)-\Gamma(a,t) \geq \frac{1}{2}(\left \lfloor \frac{t}{a}\right \rfloor^2+\left \lfloor \frac{t}{a}\right \rfloor)= 
\bin{\left \lfloor \frac{t}{a}\right \rfloor+1}{2}\geq 0$, since $r,l,s'\geq 0$ and $u^2-u\geq 0$ for all $u\geq 0$. 
Moreover, $\Gamma(a-1,t)-\Gamma(a,t) = \bin{\left \lfloor \frac{t}{a}\right \rfloor+1}{2}$ 
if and only if 
\[ 
(l,s')= 
 \left\{ 
  \begin{array}{l} 
   (0,s')\\ 
   (1,0)   
  \end{array} 
 \right. 
\] 

These two possibilities imply that $r+s<a$ and $r+s=a$, respectively. 
Finally, it is clear that $\bin{\left \lfloor \frac{t}{a}\right \rfloor+1}{2}=0$ if and only if $0\leq t < a$. 
   
{\bf (iii)} 
Let $a+t-1=a r+s$ with %$r \geq 1$ and 
$0\leq s< a$, 
then 
\[ 
a+t= 
\begin{cases}
   a r+(s+1) & \mbox{ if } 0\leq s<a-1,\\ 
   a (r+1) & \mbox{ if } s=a-1. 
\end{cases}
\] 
Hence %and by (i) we have that
\[ 
\begin{array}{cl} 
\Gamma(a,t)-\Gamma(a,t-1) & = 
 \left\{ 
  \begin{array}{l} 
(a-s-1)\bin{r}{2}+(s+1)\bin{r+1}{2}-(a-s)\bin{r}{2}-s\bin{r+1}{2}\\ 
 \\ 
a\bin{r+1}{2}-\bin{r}{2}-(a-1)\bin{r+1}{2} 
  \end{array} 
 \right.\\ 
\\   
& = \bin{r+1}{2}-\bin{r}{2} = r =\left \lfloor \frac{a+t-1}{a}\right \rfloor. 
\end{array} 
\] 
{\bf (iv)} Let $a=a_1+a_2$ and $t=t_1+t_2$, then
by $(i)$ 
\[
\begin{array}{cl} 
\Gamma(a, t) &={\rm min} \left\{ \sum_{i=1}^{a}\bin{z_i}{2}\, : \, \sum_{j=1}^a z_j= a+t \mbox{ and } z_j \geq 0 \mbox{ for all }  1 \leq j \leq a \right\}\\
& \leq (a_1-s_1)  \bin{r_1}{2}+s_1  \bin{r_1+1}{2}+ (a_2-s_2)  \bin{r_2}{2}+s_2  \bin{r_2+1}{2}= \Gamma(a_1, t_1)+\Gamma(a_2, t_2),
\end{array} 
\]
where $a_i+t_i=r_ia_i+s_i$ with $0\leq s_i< a_i$ for all $i=1,2$.

In order to have the equality in the previous inequality we need that either 
$r_1=r_2$,  $r_1=r_2+1$ and $s_1=0$ or $r_2=r_1+1$ and $s_2=0$.
%Follows directly from the definition of $\Gamma(a,t)$. 

{\bf (v}) 
Let $a+t=a r+s$ with $r \geq 1$ and $0\leq s< a$.
Thus
\[ 
\begin{array}{cl} 
\left\lceil \frac{2(a-1+\Gamma(a,t))}{a+t} \right\rceil  & =
\left\lceil \frac{2\left(a-1+(a-s)\bin{r}{2}+s\bin{r+1}{2}\right)}{a+t} \right\rceil = 
%\left\lceil \frac{2(a-1)+(a-s)r(r-1)+s(r+1)r}{a+t} \right\rceil=
\left\lceil \frac{2(a-1)+r(a r+s)-r(a - s)}{a r+s} \right\rceil\\
& \\ 
& = r+\left\lceil \frac{2(a-1)-r(a - s)}{a r+s} \right\rceil =
1+\left\lfloor \frac{t}{a} \right\rfloor +\left\lceil \frac{2(a-1)-r(a - s)}{a r+s} \right\rceil =1+\left\lfloor \frac{t}{a} \right\rfloor +L. 
\end{array} 
\] 
%Let $L=\left\lceil \frac{2(a-1)-r(a - s)}{a r+s} \right\rceil$.
Since  $a,r\geq 1$ and $0 \leq s<a$, then  $-1\leq L \leq 1$ because
\begin{center}
$2 < (a+s)(2+r)  \Leftrightarrow  -2(a r+s)< r(s-a)+2(a-1) \Leftrightarrow -1 \leq L$, and 
$2a+rs\leq2ar+s+2 \Leftrightarrow  2(a-1)-r(a-s)\leq a r+s \Leftrightarrow L \leq 1$.
\end{center}
%$2a+rs<2a+ra<3ar+2(s+1) \Leftrightarrow  2(a-1)-r(a-s)<2(a r+s) \Leftrightarrow L \leq 2$.
Moreover, 
\[
L=-1 \Leftrightarrow 2(a-1)-r(a-s)\leq -(ar+s) \Leftrightarrow 2a +s(r+1) \leq 2 \Leftrightarrow a=1 \mbox{ and } s=0.\vspace{-4mm}
\] 
%and $L=1 \Leftrightarrow 2(a-1)-r(s-a)\leq -(ar+s) \Leftrightarrow 2a +s(r+1) \leq 2 \Leftrightarrow a=1 \mbox{ and } s=0$. 
%since $L \geq 0  \Leftrightarrow -(a r+s)< r(s-a)+2(a-1) \Leftrightarrow  2 < s(r+1)+2a \Leftrightarrow s>0 \mbox{ or } s=0 \mbox{ and } a \geq 2$ and
%$L \geq 1 \Leftrightarrow 0< r(s-a)+2(a-1) \Leftrightarrow  0 < (a-1)(2-r)+r(s-1) \Leftrightarrow r=1, s\geq 1 \mbox{ or } r=2, s\geq 2$, then
%\[
%L=
%  \begin{cases}
%  -1 &\mbox{ if } a=1,\\ 
%  1 &\mbox{ if either } r=1 \mbox{ and } s\geq 1 \mbox{ or } r=2 \mbox{ and } s\geq 2,\\
%  0 &\mbox{ otherwise. }
%  \end{cases}
%\]
%\vspace{-7mm}
%Finally
%\[
%L=\left\lceil \frac{2(a-1)-r(a - s)}{a r+s} \right\rceil =
%  \begin{cases}
%  -1 &\mbox{ if } a=1,\\ 
%  %0 &\mbox{ if } a \hspace{-1.5mm}\not | \, t \mbox{ or } a|t, \mbox{ and } a \geq 2,\\
%  1 &\mbox{ if either } r=1 \mbox{ and } s\geq 1 \mbox{ or } r=2 \mbox{ and } s\geq 2,\\
%  0 &\mbox{ otherwise. }
%  \end{cases}
%\]  
%Since 
%\vspace{-3mm}
%\begin{itemize}
%%\[ 
%%\begin{array}{cl} 
%\item $L \geq -1 \Leftrightarrow -2(a r+s)< r(s-a)+2(a-1) \Leftrightarrow  2 < (a+s)(2+r) \Leftrightarrow a,r\geq 1$,
%\item $L \geq 0  \Leftrightarrow -(a r+s)< r(s-a)+2(a-1) \Leftrightarrow  2 < s(r+1)+2a \Leftrightarrow s>0 \mbox{ or } s=0 \mbox{ and } a \geq 2$,
%\item $L \geq 1 \Leftrightarrow 0< r(s-a)+2(a-1) \Leftrightarrow  0 < (a-1)(2-r)+r(s-1) \Leftrightarrow r=1, s\geq 1 \mbox{ or } r=2, s\geq 2$.
%%\end{array} 
%%\] 
%\end{itemize}
\QED 

%%%%%%%%%%%%%%%%%%%%%%%%%%%%%%%%%%%%%%%%%%%%%%%%%%%%%%%%%%%%%%%%%%%%%%%%%%%%%%%%%%%%%%%
%%%%%%%%%%%%%%%%%%%%%%%%%%%%%%%%Teorema 1-Minimal%%%%%%%%%%%%%%%%%%%%%%%%%%%%%%%%%%%%%%
%%%%%%%%%%%%%%%%%%%%%%%%%%%%%%%%%%%%%%%%%%%%%%%%%%%%%%%%%%%%%%%%%%%%%%%%%%%%%%%%%%%%%%%

\begin{Theorem}{\rm (\cite[Theorem 3.3]{BoundsIIArxiv})}\label{minimal} 
Let $G$ be a graph, then 
\[ 
q(G) \geq \alpha(G)-c(G)+\Gamma(\alpha(G), \tau(G)). 
\] 
\end{Theorem} 
\demo 
We will use induction on $\tau(G)$. 
The stars ${\cal K}_{1,n}$ ($\alpha({\cal K}_{1,n})=n-1$)  are the unique connected graphs with $\tau(G)=1$.
Since
\[ 
q({\cal K}_{1,n})=n-1=(n-1)-1+1=\alpha({\cal K}_{1,n})-c({\cal K}_{1,n})+\Gamma(n-1,1), 
\]
then the result  clearly follows.
Moreover, the stars ${\cal K}_{1,n}$ are $q$-minimal graphs. 

So we can assume that the result is true for $\tau(G)\leq k$ and $k > 1$. 
Let $G$ be a $q$-minimal graph with $\tau(G)=k+1$. 
Now, we will use induction on $\alpha(G)$. 
If $\alpha(G)=1$, then $G$ is a complete graph ${\cal K}_n$ ($\tau({\cal K}_n)=n-1$). 
Since, 
$$ 
q({\cal K}_n)=\bin{n}{2}=1-1+\bin{n}{2}=\alpha({\cal K}_n)-c({\cal K}_n)+\Gamma(1,n-1), 
$$ 
it follows that all the complete graphs satisfy the result. 
  
Hence we can assume that $\alpha(G)\geq 2$.  
Furthermore, by Lemma~\ref{lemaG}(iv),
$q(G)=\sum_{i=1}^s q(G_i)$, $\alpha(G)=\sum_{i=1}^s \alpha(G_i)$ 
and $\tau(G)=\sum_{i=1}^s \tau(G_i)$ where $G_1,\ldots,G_s$ are the 
connected components of $G$. 
Then by the induction hypothesis we can assume that $G$ is a connected graph.

%it follows from Lemma~\ref{lemaG}(iv) that we can 
%assume with out loss of generality that $G$ is connected and $\alpha(G)\geq 2$. 

Let $e$ be an edge of $G$ and consider the graph 
$G'=G-e$. We have two possibilities 
\[ 
\tau(G')= 
 \left\{ 
  \begin{array}{l} 
   \tau(G)\\ 
   \tau(G)-1 
  \end{array} 
 \right. 
\] 
That is, an edge of $G$ is either a bridge or critical. 
  
\medskip 

\noindent {\bf Case 1} First, assume that $G$ has no bridges, that is, $G$ is an $\alpha$-critical graph. 
\begin{Claim}\label{claim0}
Let $v$ be a vertex of $G$ of maximum degree, then
\[
\deg(v')\geq 1+ \left\lfloor \frac{\tau(G)-1}{\alpha(G)}\right\rfloor. 
\] 
\end{Claim}
\demo
Since any $\alpha$-critical graph is $\tau$-critical, then $\tau(G-v)=\tau(G)-1$ and $\alpha(G-v)=\alpha(G)$.
Moreover, since the $\alpha$-critical graphs are $2$-connected, then $G-v$ is connected. 
Now, by the induction hypothesis we have that 
\[ 
q(G-v)\geq \alpha(G)-1+\Gamma(\alpha(G),\tau(G)-1).
\] 
Using the formula 
$$ 
\sum_{v_i \in V(G-v)} \deg(v_i)=2q(G-v) 
$$
and Lemma~\ref{lemaG}(v), we conclude that there must exist a vertex $v'\in V(G-v)$ with 
\begin{eqnarray}\label{eq10} 
\deg(v')\geq \left\lceil \frac{2q(G-v)}{|V(G-v)|}\right\rceil 
&\geq& \left\lceil \frac{2(\alpha(G)-1+\Gamma(\alpha(G),\tau(G)-1))}{n-1}\right\rceil \nonumber\\
&\overset{\rm (v)}{\geq}& 1+ \left\lfloor \frac{\tau(G)-1}{\alpha(G)}\right\rfloor+L. 
\end{eqnarray} 

Now, since $L\geq 0$ ($\alpha(G)\geq 2$), then by Lemma~\ref{lemaG}(iii) %we obtain that  
\begin{eqnarray}\label{eq4} 
q(G) =  q(G-v)+\deg(v) &\geq&  \alpha(G)-1+\Gamma(\alpha(G),\tau(G)-1)+\deg(v')\nonumber\\
                        & \overset{\scriptsize \rm (\ref{eq10})}{\geq} & \alpha(G)-1+\Gamma(\alpha(G),\tau(G)-1)+ 1+\left\lfloor \frac{\tau(G)-1}{\alpha(G)}\right\rfloor \nonumber\\ 
                        &\overset{\rm (iii)}{=}& \alpha(G)-1+\Gamma(\alpha(G),\tau(G)).  
\end{eqnarray} 
%Since $L\geq 0$ (), using Lemma~\ref{lemaG}(v)
%\[
%q(G)\geq  \alpha(G)-1+\Gamma(\alpha(G),\tau(G)-1)-1
%\]

So, if the graph $G$ has an edge that is a bridge, then $c(G')=c(G)+1=2$. 
Let denote by $G_1$ and $G_2$ the connected components of $G-e$. 
Now, we need to considerer another two cases: %two more cases:
  
\noindent {\bf Case 2} Assume that $\tau(G_1)>0$ or $\tau(G_2)>0$, then 
$\tau(G_1)\leq k$, $\tau(G_2)\leq k$, and by the induction hypothesis 
\[ 
\begin{array}{c} 
q(G_1)  \geq  \alpha(G_1)-1+\Gamma(\alpha(G_1),\tau(G_1))\mbox{ and }
q(G_2)  \geq  \alpha(G_2)-1+\Gamma(\alpha(G_2),\tau(G_2)).
\end{array} 
\] 
Using the above formulas, the fact that $\alpha(G)=\alpha(G_1)+\alpha(G_2)$ 
and $\tau(G)=\tau(G_1)+\tau(G_2)$, and Lemma~\ref{lemaG}(iv), we get 
\[ 
\begin{array}{cl} 
q(G) & = q(G_1)+q(G_2)+1\\ 
& \geq  \alpha(G_1)-1+\alpha(G_2) -1+\Gamma(\alpha(G_1),\tau(G_1))+\Gamma(\alpha(G_2),\tau(G_2))+1\\ 
& = \alpha(G)-1 +\Gamma(\alpha(G_1),\tau(G_1))+\Gamma(\alpha(G_2),\tau(G_2)) 
\overset{\rm (iv)}{\geq} \alpha(G)-1+\Gamma(\alpha(G),\tau(G)).
\end{array} 
\] 
\newpage 

\noindent {\bf Case 3} Assume that there does not exist a bridge satisfying the 
above conditions: for all the bridges of $G$ we have that $\tau(G_1)=0$ or $\tau(G_2)=0$. 
That is, each bridge connects an isolated vertex with the rest of $G$.
In this case, we must have that $G$ is 
equal to an $\alpha$-critical graph $G_1$ with a vertex of $G_1$ being the center 
of a star ${\cal K}_{1,l}$. 
Moreover, $\tau(G)=\tau(G_1)$ and $\alpha(G)=l+\alpha(G_1)$ because $G_1$ is vertex-critical and 
therefore each vertex belongs to a minimum vertex cover. 
Now using Case 1 and Lemma~\ref{lemaG}(ii), we obtain, 
\[ 
\begin{array}{cl} 
q(G)=l+q(G_1) & \geq l+\alpha(G_1)-1+\Gamma(\alpha(G_1),\tau(G_1))\\ 
& = \alpha(G)-1+\Gamma(\alpha(G_1),\tau(G))
 \overset{\rm (ii)}{\geq} \alpha(G)-1+\Gamma(\alpha(G),\tau(G)). 
\end{array} \vspace{-7mm}
\] \QED

\paragraph{The $k$-connexion of a graph}
A {\it $k$-connexion} of a graph $G$ is a $k$-connected graph $G'$
on the same vertex set as $G$, with the minimum possible number of edges, 
and such that $G$ is a subgraph of $G'$.
The graph $G$ is called the {\it subjacent} graph of the $k$-connexion graph $G'$ and the
edges of $G'$ that are not edges of $G$ are called the {\it connexion edges}.

Clearly a $1$-connexion graph $G'$ of a disconnected graph $G$ can be obtained by adding $c(G)-1$
edges, where $c(G)$ is the number of connected components of $G$.
This definition is equivalent 
to the one given in~\cite{bougard} of a {\it tree-linking} of a graph. 
In fact, $2$-connexions as defined here, are equivalent to 
{\it cycle-linkings} as defined in~\cite{bougard}. 

\begin{Example}\rm
In order to illustrate the concept of $1$-connexion consider the following graphs:

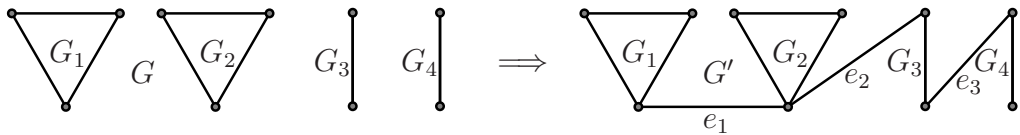
\vspace{13mm} 
\begin{figure}[h] 
\begin{center} 
$ 
\setlength{\unitlength}{.40mm} 
\thicklines\begin{picture}(335,0) 

\scalebox{.47}{\includegraphics{PTT1e}} 
{\large 
\put(17,14){$\Longrightarrow$} 
\put(-132,17){$G_1$} 
\put(-82,17){$G_2$} 
\put(-105,10){$G$} 
\put(-44,14){$G_3$} 
\put(-15,14){$G_4$} 
} 

\hspace{15mm}
\scalebox{.47}{\includegraphics{1conexione}}
{\large 
\put(-132,17){$G_1$} 
\put(-82,17){$G_2$} 
\put(-105,10){$G'$} 
\put(-44,14){$G_3$} 
\put(-15,14){$G_4$}
\put(-105,-5){$e_1$}
\put(-58,10){$e_2$}
\put(-21,8){$e_3$}
}   
\end{picture} 
$ 
\end{center} 
\end{figure}

\vspace{-5mm}

\noindent $G'$ is the $1$-connexion of $G$, the edges $e_1,e_2,e_3$ are the $1$-connexion edges of $G'$, 
and $G_1,G_2,G_4$ are the leaves of $G'$.
\end{Example}

A {\it leaf} of a $1$-connexion $G'$ of a graph $G$ is a connected component $G_i$ of $G$ 
with the property that there exists a unique vertex $v$ of $G_i$, such that all connexion 
edges with one end in $G_i$ are incident to $v$.
If $G$ is a connected graph, then we say that $G$ is a {\it leaf} of $G'$.
Note that a $1$-connexion $G'$ of a graph $G$ has at least one leaf. 
%either a connected component $G_i$ of $G$ incident to a unique connexion edge 
%or a connected component $G_i$ with the property that there exists a unique vertex 
%$v$ in $G_i$ such that all connexion edges with one end in $G_i$ are incident to the vertex $v$.
\paragraph{Polygon transformed Tur\'an graph}
A graph $G$ with covering number $\tau(G)=t$ and stability number $\alpha(G)=a$ is said to be a 
{\it polygon transformed Tur\'an graph} or {\it PTT graph} if either $G$ is 
isomorphic to $T(a, t)$, or $a \leq t < 2a$ and $G$ can be obtained from 
$T(a, t)$ by the following construction:
 
Let $k_2$ and $k_3$ be the number of copies of $K_2$ and $K_3$  in $T(a, t)$ respectively.
%For all $i$ such that $i \leq \min\{k_2,k_3\}$ we replace $j_i$ copies of $K_2$ 
%and one copy of $K_3$ by a cycle $C_{2j_i+3}$,
%where $j_1+\cdots+j_k \leq k_2$.
Let $k$ be a positive integer with $k \leq \min\{k_2,k_3\}$ and 
for all $1\leq i\leq k$ take positive integers $j_i$ such that $j_1+\cdots+j_k \leq k_2$. 
%where $k_2$,$k_3$ denote the number of copies in $T(a, t)$ 
%of $K_2$ and $K_3$ respectively.
Finally, for all  $i=1,\ldots,k$ 
replace one copy of $K_3$ and $j_i$ copies of $K_2$ by a cycle $C_{2j_i+3}$.
%for all  $i=1,\ldots,k$.

\medskip

In this way a PTT graph is the disjoint union of complete graphs and possibly odd cycles.

%\newpage
\begin{Example}\rm
In order to illustrate the previous concept consider the following graphs:
\vspace{16mm} 
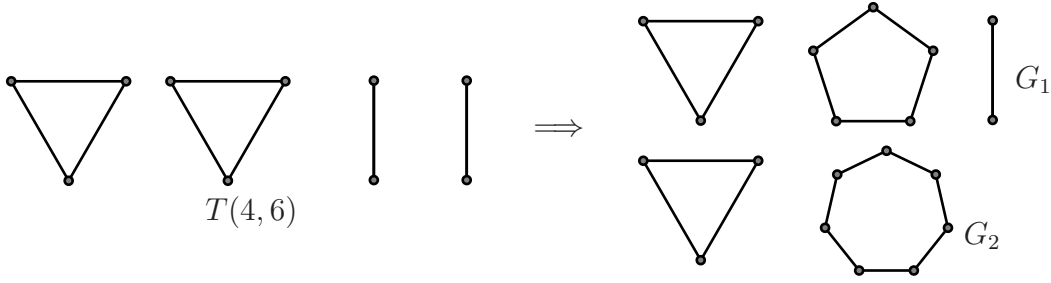
\begin{center} 
$ 
\setlength{\unitlength}{.40mm} 
\thicklines\begin{picture}(335,0) 
\scalebox{.50}{\includegraphics{PTT1e}} 
{\large 
\put(20,17){$\Longrightarrow$} 
\put(-89,-11){$T(4,6)$} 
\put(180,33){$G_1$} 
\put(163,-19){$G_2$} 
} 
\end{picture} 
$ 
\end{center} 

\vspace{-19mm} 
\begin{center} 
$ 
\setlength{\unitlength}{.40mm} 
\thicklines\begin{picture}(-85,0) 
\scalebox{.50}{\includegraphics{PTT2e}}  
\end{picture} 
$ 
\end{center} 

\vspace{8mm} 
\begin{center} 
$ 
\setlength{\unitlength}{.40mm} 
\thicklines\begin{picture}(-85,0) 
\scalebox{.50}{\includegraphics{PTT3e}}  
\end{picture} 
$ 
\end{center} 

\vspace{-2mm}
\noindent 
in the left side it can be seen the Tur\'an graph $T(4,6)$ and in the right side there are two of the three possible polygon transformed graph of $T(4,6)$.
Note that $G_1$ and $G_2$ are obtained when we take $1=k<\min\{k_2,k_3\}=2$, and $j_1=1$ and $j_1=2$ respectively.
\end{Example}

%%%%%%%%%%%%%%%%%%%%%%%%%%%%%%%%%%%%%%%%%%%%%%%%%%%%%%%%%%%%%%%%%%%%%%%%%%%%%%%%%%%%%%%
%%%%%%%%%%%%%%%%%%%%%%%%%%%%%Clasificacion 1-Conexas%%%%%%%%%%%%%%%%%%%%%%%%%%%%%%%%%%%
%%%%%%%%%%%%%%%%%%%%%%%%%%%%%%%%%%%%%%%%%%%%%%%%%%%%%%%%%%%%%%%%%%%%%%%%%%%%%%%%%%%%%%%

\begin{Lemma}\label{conexo}
A graph $G$ is $q$-minimal if and only if $G$ is a $1$-connexion 
of a polygon transformed Tur\'an graph.
\end{Lemma}
\demo
$(\Leftarrow)$
Let $H$ be a PTT graph with $H_1,\ldots,H_a$ connected components and let $L$ be a $1$-connexion of $H$.
Since ${\cal K}_r$ and $C_{2s+1}$ are $q$-minimal, $\alpha(L)=\sum_{i=1}^a \alpha(H_i)$, $\alpha(C_{2s+1})=s$, 
%$q({\cal K}_r)=\alpha({\cal K}_r)-1+\Gamma(\alpha({\cal K}_r), \tau({\cal K}_r))$, 
%$q(C_{2s+1})=2s+1=\alpha(C_{2s+1})-1+\Gamma(\alpha(C_{2s+1}), \tau(C_{2s+1}))$, and 
$q(L)=a-1+\sum_{i=1}^a q(H_i)$, then $L$ is $q$-minimal.

$(\Rightarrow)$ 
We  use double induction on the stability and covering numbers of the graph.
For $\alpha(G)=1$, $G$ must be a complete graph and the result is clear.

Let $G$ be a $q$-minimal graph with $\alpha(G) \geq 2$. 
We can assume that $G$ is an $\alpha$-critical graph, since
if $G$ is not $\alpha$-critical, then using the arguments used in cases 2 and 3 
(in the proof of Theorem~\ref{minimal}) and the induction hypothesis, the result follows readily. 
%Therefore we can assume that $G$ is $\alpha$-critical.
Since %a $1$-connexion of a PTT graph is a $q$-minimal graph, 
%then the 
connexion edges of a $1$-connexion of a disconnected PTT graph are not $\alpha$-critical edges
($\alpha(L)=\sum_{i=1}^k \alpha(H_i)$, where $L$ is a $1$-connexion of a disconnected PTT graph $H$ 
with connected components $H_1,\ldots, H_k$), 
%Therefore 
then the result follows if we prove that $G$ is either an odd cycle or a complete graph.

\begin{Claim}\label{claim1} 
If %$G$ be a $q$-minimal and $\alpha$-critical graph and let 
$v$ is a vertex of $G$ of maximal degree, then $G-v$ is $q$-minimal.
\end{Claim}
\demo
Assume that $G-v$ is not $q$-minimal.
Then, by Claim~\ref{claim0} %the same arguments as those in case 1 (in the proof of Theorem~\ref{minimal}) 
and Lemma~\ref{lemaG},
\begin{eqnarray}\label{eq5}
q(G) = q(G-v)+\deg(v) %&\geq & q(G-v)+\deg(v') \nonumber \\
                & \overset{\rm (v)}{\geq} & 
                \alpha(G)+\Gamma(\alpha(G),\tau(G)-1)+\left \lfloor \frac{\alpha(G)+\tau(G)-1}{\alpha(G)} \right \rfloor \nonumber\\
                & \overset{\rm (iii)}{=} & \alpha(G)+\Gamma(\alpha(G),\tau(G)); 
\end{eqnarray}
which is a contradiction to the $q$-minimality of $G$. 
\QED

Now, since $G-v$ is $q$-minimal, then by the induction hypothesis $G-v$ is a $1$-connexion of a PTT graph.
Moreover, since $G$ is $\alpha$-critical, then $G\setminus N[v]$  (where $N[v]=N(v)\cup \{v\}$) is a maximal induced subgraph of $G$
with $\alpha(G\setminus N[v])=\alpha(G)-1$. %for all $v\in V(G)$.
%(in particular $\tau$-critical), then $\alpha(G)=\alpha(G-v)$, 
%that is, the set of vertices $N(v)$ must satisfy that $N(v)\cap M\neq \emptyset$ 
%for all maximum stable sets $M$ of $G-v$.
%Hence, $\alpha(G\setminus N[v])=\alpha(G)-1$, where $N[v]=N(v)\cup \{v\}$. 
%A set of vertices $N$ in $G-v$ can be the set of neighbors of $v$ in $G$   
%if and only if $V(G-v) \setminus N$ induces a subgraph $G'$ of 
%$G-v$ with $\alpha(G')=\alpha(G)-1$. 
%Moreover, $N$ is minimal under inclusion if and only if $G[V(G-v) \setminus N]$ is maximal under inclusion.
Therefore, we need to determine the maximal induced subgraphs $L'$ of a $1$-connexion of a PTT graph $L$ with $\alpha(L')=\alpha(L)-1$.
%In this case we can determine the maximal induced subgraphs $G'$ of $G-v$ with $\alpha(G')=\alpha(G)-1$.

\begin{Claim}\label{claim2}
Let $H$ be a PTT graph with $H_1,\ldots,H_a$ connected components and $L$ be a $1$-connexion of $H$. 
If $L'$ is a maximal induced subgraph of $L$ with $\alpha(L')=\alpha(L)-1$, then
\begin{description}
\item[{\it (i)}] $L'$ is induced by the set of vertices $V(L)\setminus V(H_i)$, for some $H_i$ with $\alpha(H_i)=1$, or

\item[{\it (ii)}] $L'$ is induced by the set of vertices in $V(L)\setminus \{v_1,v_2,v_3\}$, 
                  where $\{v_1,v_2,v_3\}$ are vertices of an odd cycle $H_j$ such that $H_j\setminus \{v_1,v_2,v_3\}$ 
                  is a disjoint union of paths with an even number of vertices, or

\item[{\it (iii)}] $L'$ satisfies the following conditions: 
                   ${\it (1)}$ $V(H_i)\cap V(L')\neq \emptyset$ for all $H_i$, 
                   ${\it (2)}$ if $H_i$ is an odd cycle, then $V(H_i)\subset V(L')$, 
%                   ${\it (3)}$ if $H_i$ is a complete graph, %such that $V(H_i) \not \subset V(L')$, 
%                  then for all $v \in V(H_i)\cap V(L')$
%                  there exists at least one connexion edge $e_v$ of $L$ incident to $v$.                                       
                   ${\it (3)}$ if $H_i$ is a complete graph such that $V(H_i) \not \subset V(L')$, then for all $v \in V(H_i)\cap V(L')$
                   there exists at least one connexion edge $e_v$ of $L$ incident to $v$.                   
\end{description}
\end{Claim}
\demo
%We divide the proof in the following cases:
%
%{\bf Case ( $V(L) \cap V(H_j)=\emptyset$ for some $H_j$ with $\alpha(H_j)=1$)}
If $V(L')\cap V(H_i)=\emptyset$ for some $1\leq i\leq a$ with $\alpha(H_i)=1$, then 
%In this case, 
$L'=L[V(L)\setminus V(H_i)]$, 
since $V(L') \subseteq V(L) \setminus V(H_i)$ and 
$\alpha(L[V(L)\setminus V(H_i)])=\alpha(L)-1$.

\vspace{17mm} 
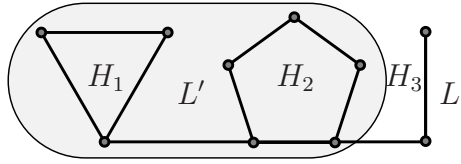
\begin{center} 
$ 
\setlength{\unitlength}{.40mm} 
\thicklines\begin{picture}(150,0) 
\scalebox{.55}{\includegraphics{PTT(i)e.eps}} 
{\large 
\put(-85,19){$L'$} 
\put(2,19){$L$} 
\put(-115,25){$H_1$} 
\put(-52,25){$H_2$} 
\put(-16,25){$H_3$} 
} 
\end{picture} 
$ 
\end{center} 
\vspace{-2mm}

Therefore we can assume that %if $L'$ does not satisfy $(i)$, and thus 
$V(H_i)\cap V(L')\neq \emptyset$ for all $H_i$ with $\alpha(H_i)=1$.
%{\bf Case ( $V(L) \cap V(H_j) \neq \emptyset$ for all $H_j$ with $\alpha(H_j)=1$)}
Now, let assume that $V(H_j) \not\subset V(L')$ for some $H_j=C_{2m+1}$. 
%Let $H_j$ be an odd cycle with $2m+1$ vertices.
Since all the proper induced graphs of a cycle 
are paths $P_{n}$ with $\alpha(P_n)=\lceil \frac{n}{2} \rceil$, 
then $\alpha(H_j\setminus C)=\alpha(H_j)-1=m-1$ for some $C\subset V(H_j)$ if and only if
$H_j[V(H_j)\setminus C]$ is a disjoint union of three paths $P_{m_1},P_{m_2},P_{m_3}$ 
for some even numbers $m_1,m_2,m_3\geq 0$ such that $m_1+m_2+m_3=2(m-1)$.
Since $L'$ is a maximal induced subgraph of $L$ with $\alpha(L')=\alpha(L)-1$,
then %$\alpha(L[V(H_j)\cap V(L')])$ is equal to either $m$ or $m-1$. Moreover, 
$V(H_j) \not\subset V(L')$ for only one $H_j=C_{2m+1}$; therefore $L'$ is given by $(ii)$.
%Since $\alpha(L\setminus C)=\alpha(L)-1$, where $H_j\setminus C$ 
%is a disjoint union of paths with an even number of vertices, then $L'=L\setminus C$.

\vspace{17mm} 
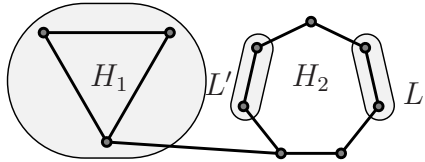
\begin{center} 
$ 
\setlength{\unitlength}{.40mm} 
\thicklines\begin{picture}(140,0) 
\scalebox{.55}{\includegraphics{PTT(ii)e}} 
{\large 
\put(-64,22){$L'$} 
\put(2,19){$L$} 
\put(-102,25){$H_1$} 
\put(-35,25){$H_2$}  
} 
\end{picture} 
$ 
\end{center} 
\vspace{-2mm}

%Therefore, either $L'$ satisfies {\it (ii)} or 
%$V(H_j)\subset V(L')$ for all $H_j$ with $\alpha(H_j)\geq 2$.
In order to finish, %if $L'$ is not given by ${\rm (i)}$ or ${\rm (ii)}$, then 
we can assume that $V(L')\cap V(H_i)\neq \emptyset$ for all $1\leq i\leq a$
and $V(H_j)\subset V(L')$ for all $H_j$ with $\alpha(H_j)\geq 2$. 
%Moreover, if $H_j$ is an odd cycle, then $V(H_j)\subset V(L')$.
Clearly, if $v\in V(H_i)\cap V(L')$ such that $v$ is not incident to any connexion edge of $L$,
then %$V(H_i) \not \subset V(L')$, then 
$V(H_i)\subset V(L')$
because $\alpha(L')=\alpha(L[V(L')\cup V(H_i)])$.
%Moreover, if $V(H_i) \subset V(L')$, then $\alpha(L')=\alpha(L)$; a contradiction.
\vspace{13mm} 
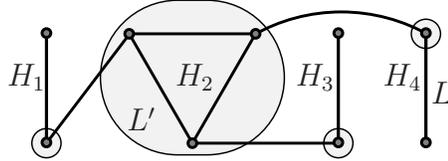
\begin{center} 
$ 
\setlength{\unitlength}{.40mm} 
\thicklines\begin{picture}(145,0) 
\scalebox{.55}{\includegraphics{PTT(iii)e}} 
{\large 
\put(-105,10){$L'$} 
\put(-4,19){$L$} 
\put(-145,25){$H_1$} 
\put(-89,25){$H_2$} 
\put(-49,25){$H_3$}
\put(-20,25){$H_4$} 
} 
\end{picture} 
$ 
\end{center} 
\vspace{-7mm}
\QED

%\begin{Example}
%baba
%
%\noindent $(i)$
%$L'=L[\{v_1,\ldots,v_8\}]$
%
%\vspace{17mm} 
%\begin{center} 
%$ 
%\setlength{\unitlength}{.40mm} 
%\thicklines\begin{picture}(335,0) 
%\scalebox{.70}{\includegraphics{PTT(i)e}} 
%{\large 
%\put(-66,17){$T(4,6)$} 
%} 
%\end{picture} 
%$ 
%\end{center} 
%\end{Example}

\vspace{6mm}

Take $L=G\setminus v$ and $L'=G\setminus N[v]$.
If $G\setminus N[v]$ satisfies ${\rm (i)}$, then $G$ must be a complete graph.
If $G\setminus N[v]$ satisfies ${\rm (ii)}$, then $G$ must be an odd subdivision of the complete graph ${\cal K}_4$
and it fails in one edge in order to be a $q$-minimal graph.
%Applying Claim~\ref{claim2} to $G-v$ it is easy to conclude that
%\begin{description}
%\item[$\bullet$] $G$ is a complete graph whenever $G\setminus N[v]$ satisfies {\it (i)}.
%\item[$\bullet$] $G$ is not $q$-minimal whenever $G\setminus N[v]$ satisfies {\it (ii)}.
%\end{description}

%Therefore it only remains to considerer the case when $G\setminus N[v]$ satisfies ${\it (iii)}$.
Finally, assume that $G\setminus N[v]$ satisfies ${\rm (iii)}$.
Let $H_{i_0}$ be a connected component of $H$ such that $H_{i_0}$ is a complete graph and
%$H_{i_0}$ is a connected component of the subjacent graph of 
%$G-v$ (a PTT graph) with 
$V(H_{i_0}) \not \subset V(G\setminus N[v])$. 
%(note that by Claim~\ref{claim2} (iii) there exists at 
%least one graph $H_i$ with this condition) and 
Take $P=V(H_{i_0}) \cap V(G\setminus N[v])$ and $Q=V(H_{i_0})\setminus P$.
Since $G-v$ is $q$-minimal, then for all $u\in P$, the graph $(G-v)\setminus u$ is disconnected.
For all $u\in P$, let $G_u$ be the disjoint union of the 
connected components $C_1,\ldots, C_s$ ($C_i$ is a $1$-connexion of some PTT graph) 
of $(G-v)\setminus u$ with $V(C_j) \cap V(H_{i_0})=\emptyset$. %for all $1\leq j\leq s$.
%That is, $C_i$ is a $1$-connexion of some PTT graph. %and in general $G_u$ is not a connected graph.
Note that, $G_u$ is an induced subgraph %(a disjoint union of a $1$-connexion of PTT graphs) 
of $G-v$ and its connected components are joined to $u$ by some connexion edges.
%Moreover, if $| V(H_{i_0})| \geq 2$, then $G_u$ is unique.

%Here we need to considerer two cases, the first case is when $G_u$ is not a PTT graph. 
Let $C$ be a connected component of $G_u$ and $S$ be a leaf of $C$ not joined to $u$ by a connexion edge. 
Since $G$ is a $2$-connected graph, then $v$ must be incident with at least one vertex of $S$.
If $G_u$ is either a complete graph or an odd cycle, then $v$ must be incident with at least one vertex of $G_u$
not incident with a connexion edge of $G-v$.
%In the other case, if $G_u$ is a PTT graph, 
%then by the $2$-connectivity of $G$, there exists at least one vertex $w$ such 
%that $w$ is incident to $v$ and we can consider that $S=G_u$ is the only leaf of $G_u$. 
Moreover, if $v_s$ is the unique vertex of a leaf $S$ of $C$ such that all the connexion edges 
with one end in $S$ are incident to $v_s$, then by Claim~\ref{claim2} (iii), 
$v$ must be incident with all the vertices of $S\setminus v_s$. 
Since $v$ is incident with all the vertices of $Q$, then
\begin{eqnarray}\label{eq6}
\deg(v) \geq |Q|+\sum_{u\in P} \sum_{ H_j \in L(G_u)} (|H_j|-1) \overset{(\ast)}{\geq} |H_{i_0}|,
\end{eqnarray}
where $L(G_u)$ is either the set of leaves of $G_u$ not joined to $u$ or if $G_u$ is a $2$-connected graph, then $L(G_u)$ is $\{G_u\}$. 
%when $G_u$ is not a PTT graph or equal to $G_u$ when $G_u$ is a PTT graph.
Furthermore, $(\ast)$ is an equality %$\deg(v)=|H_{i_0}|$ 
if and only if $G_u$ is connected, all the leaves of $G_u$ are isomorphic to 
${\cal K}_2$, and if $G_u$ has at most two leaves. %and if  $G_u$ is not a PTT graph, then $G_u$ has exactly two leaves.

\vspace{37mm} 
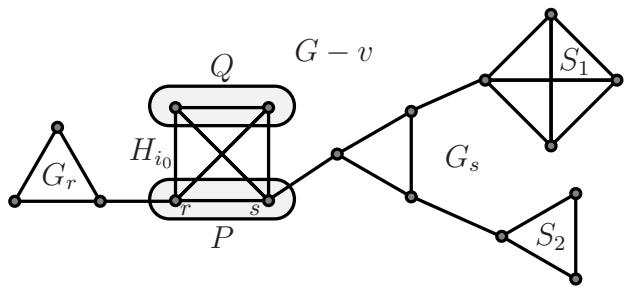
\begin{center} 
$ 
\setlength{\unitlength}{.40mm} 
\thicklines\begin{picture}(205,0) 

\scalebox{.62}{\includegraphics{figura5e}} 

{\footnotesize
\put(-148,24){$r$} 
\put(-125,24){$s$} 
}
{\large 
\put(-194,34){$G_r$} 
\put(-60,40){$G_s$} 
\put(-165,42){$H_{i_0}$}
\put(-22,72){$S_1$} 
\put(-30,14){$S_2$} 
\put(-138,13){$P$} 
\put(-138,70){$Q$} 
\put(-110,75){$G-v$} 
} 
\end{picture} 
$ 
\end{center} 
\vspace{-3mm} 

%\medskip 

Using the first inequality in the equation~(\ref{eq6}), it is not difficult to prove that
%\[\deg(v)\geq s+(k-s)(k-1)\geq 2k-2, \mbox{ where }k=\lfloor \frac{|V(G-v)|}{\alpha(G-v)}\rfloor.\]
\begin{equation}\label{eq9}
\deg(v)\geq s+(k-s)(k-1)\geq 2k-2, \mbox{ where }k=\left \lfloor \frac{|V(G-v)|}{\alpha(G)} \right \rfloor \leq |V(H_{i_0})|.
\end{equation}

On the other hand, since $G$ and $G-v$ are $q$-minimal graphs, then
%Now, let $H_{i_0}$ with $|H_{i_0}|=k=\lfloor \frac{|V(G-v)|}{\alpha(G)}\rfloor$, using that
%\vspace{-1mm} 
%\[
$
\deg(v)=q(G)-q(G-v) = \Gamma(\alpha(G),\tau(G))-\Gamma(\alpha(G),\tau(G)-1)
= \left\lfloor \frac{|V(G-v)|}{\alpha(G)} \right \rfloor= \left \lceil \frac{|V(G)|}{\alpha(G)} \right \rceil -1\leq k.
$
%\]

%Now, if $|H_{i_0}|=\lfloor \frac{|V(G-v)|}{\alpha(G)}\rfloor=k$, 
%then by Equation~(\ref{eq6}), $\deg(v)\geq 2k-2$.
Therefore, $k=2$,
%A similar argument shows that $k=2$ when we take $|H_{i_0}|=k=\lceil \frac{|V(G-v)|}{\alpha(G)}\rceil$.
$H_{i_0}={\cal K}_2$, $\deg(v)=2$, $ \left\lceil \frac{|V(G-v)|}{\alpha(G)} \right \rceil\leq 3$, 
$|P|=1$, and $G_u$ has only two leaves, 
that is, $G$ is an odd cycle. %because $G_u$ is a $1$-connexion of a 
%PTT graph whose components are all isomorphic to ${\cal K}_2$.
\QED
%%%%%%%%%%%%%%%%%%%%%%%%%%%%%%%%%%%%%%%%%%%%%%%%%%%%%%%%%%%%%%%%%%%%%%%%%%%%%%%%%%%%%%%%%%%%%%%%%%
%%%%%%%%%%%%%%%%%%%%%%%%%%%%%%%%%%%%%%%%2-connected case%%%%%%%%%%%%%%%%%%%%%%%%%%%%%%%%%%%%%%%%%%
%%%%%%%%%%%%%%%%%%%%%%%%%%%%%%%%%%%%%%%%%%%%%%%%%%%%%%%%%%%%%%%%%%%%%%%%%%%%%%%%%%%%%%%%%%%%%%%%%%

\subsection{The $2$-connected case}

\begin{Theorem}\label{2conexo}
Let $G$ be a $2$-connected graph with $\tau(G)\geq 2$, then 
\[
q(G) \geq 
   \begin{cases}
      2\alpha(G) & \mbox{ if } \tau(G) \leq \alpha(G), \\
      \alpha(G)-1+\Gamma(\alpha(G), \tau(G))  & \mbox{ if } \alpha(G)=1 \mbox{ or } \tau(G)-\alpha(G)=1,\\ 
      \alpha(G)+\Gamma(\alpha(G), \tau(G)) & \mbox{ otherwise.} 
   \end{cases}
\]
Furthermore, $G$ is a $\{q,2\}$-minimal with $\tau(G)\geq \alpha(G)$ if and only if one of the following conditions is satisfied:
\begin{description}
\item[{\it (i)}] $G$ is an even cycle, 
\item[{\it (ii)}] $G$ is the complete graph with at least three vertices,
\item[{\it (iii)}] $G$ is an odd cycle, 

\item[{\it (iv)}] $G$ is a $2$-connexion of a polygon transformed Tur\'an graph,

\item[{\it (v)}] $G$ is an odd subdivision of the complete graph ${\cal K}_4$,

\item[{\it (vi)}] $G$ is isomorphic to the following graph:
\vspace{21mm} 
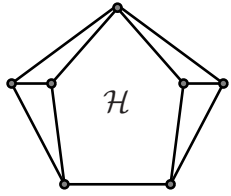
\begin{center} 
$ 
\setlength{\unitlength}{.36mm} 
\thicklines\begin{picture}(100,0) 

\scalebox{.50}{\includegraphics{c5rr}} 

{\large 
\put(-47,29){${\cal H}$}  
} 
\end{picture} 
$ 
\end{center} 
\vspace{-4mm}                     

\end{description} 
\end{Theorem}

\demo %$(\Leftarrow)$ Is not difficult to see that all graphs enumerated in the theorem are $2$-connected and satisfies the bound on the number of edges with equality.
%Therefore all these graphs are $\{q,2\}$-minimal.
%$(\Rightarrow)$
Let $H$ be a $2$-connected graph and $G$ be a $\{q,2\}$-minimal graph.
We will divide the proof in three cases: 
%we begin with the most simple cases. 

{\bf Case 1 ($\tau(G)=\alpha(G)>1$).}
Let $H$ %be a $2$-connected graph 
with $\tau(H)=\alpha(H)>1$, then ${\rm deg}(v)\geq 2$ for all $v\in V(H)$.
Therefore
\[%2\alpha(H)\geq 
q(H) = \frac{\sum_{v\in V(H)} {\rm deg}(v)}{2}\geq |V(H)|=2\alpha(H).\]

Since the even cycle $C_{2a}$ is a $2$-connected graph with 
$\alpha(C_{2a})=\tau(C_{2a})$, then $2\alpha(G) \leq q(G)\leq 2\alpha(G)$ %for all the %$\{q,2\}$-minimal graphs with $\tau(H)=\alpha(H)>1$.
and ${\rm deg}(v)=2$ for all $v\in V(G)$. %and $q(G)=2\alpha(G)$.
Furthermore, since a graph $H$ with all its vertices of degree two is a disjoint union of cycles, then
$G$ is a $\{q,2\}$-minimal graph with $\alpha(G)=\tau(G)$ if and only if $G$ is an even cycle. 

{\bf Case 2 ($\tau(G) < \alpha(G) $ and $\alpha(G) >1$).}
Since ${\rm deg}(v)\geq 2$ for all $v\in V(H)$ ($H$ is $2$-connected), 
then $q(H)\geq 2\alpha(H)$. %where $M$ is a maximum stable set of $G$. %2|M|=
Let  $s({\cal K}_{a,t})$ be an odd subdivision of the complete bipartite graph ${\cal K}_{a,t}$.
Since $s({\cal K}_{a,t})$ is a $2$-connected graph with $\alpha(s({\cal K}_{a,t}))=\tau(s({\cal K}_{a,t}))+(a-t)$ and $q({\cal K}_{a,t})=2\alpha({\cal K}_{a,t})$, 
then $q(G)= 2\alpha(G)$ for all the $\{q,2\}$-minimal graphs $G$ with $\tau(G) < \alpha(G) >1$. %for all $\{q,2\}$-minimal graph with $\tau(G)<\alpha(G)>1$.
%Moreover,  
%That is, $q(G)=2\alpha(G)$ for all the $\{q,2\}$-minimal graphs with $\tau(G)<\alpha(G)>1$.

This finishes the proof of the lower bound for the number of edges whenever $\tau(G)\leq \alpha(G)$ and $\alpha(G) >1$
and the characterization
of the $\{q,2\}$-minimal graphs (case {\rm (i)}) %, (ii)} and {\it (iii)})  
whenever  %either $\alpha(G)=1$ or 
$\alpha(G)= \tau(G)$ and $\alpha(G) >1$. %\leq \alpha(G)+1$.

\medskip

{\bf Case 3 ($\alpha(G)=1$ or $\tau(G)\geq \alpha(G)+1$).}
We  use double induction on the stability and covering numbers of the graph. %in the last two cases.
%For $\alpha(G)=1$, $G$ must be a complete graph and the result is clear.

%{\bf Case 3 ($\alpha(G)=1$ or $\tau(G)= \alpha(G)+1$).}
If $\alpha(G)=1$, then $G$ is a complete graph and clearly $G$ is a $\{q,2\}$-minimal graph.
If $\tau(H)-\alpha(H)=1$, then by Lemma~\ref{conexo}, $q(H) \geq \alpha(H)-1+\Gamma(\alpha(H), \tau(H))$ 
and $H$ is a $q$-minimal graph with $\tau(H)-\alpha(H)=1$ if and only if $H$ is an odd cycle.
Since the odd cycles are $2$-connected, %then if $H$ is a $2$-conected graph with $\tau(H)-\alpha(H)=1$,
then $q(H) \geq \alpha(H)-1+\Gamma(\alpha(H), \tau(H))$ whenever $\tau(H)= \alpha(H)+1$ and $G$ is $\{q,2\}$-minimal with $\tau(G)= \alpha(G)+1$ if and only if $G$ is an odd cycle.
%By Lemma~\ref{conexo}, $q(G)=\alpha(G)-1+\Gamma(\alpha(G), \tau(G))$ whenever 
%$G$ is a graph with $\alpha(G)=1$ or $\tau(G)-\alpha(G)=1$. 
%Moreover, Lemma~\ref{conexo} proves that a $\{q,2\}$-minimal graph satisfies these conditions 
%if and only if $G$ is a complete graphs or an odd cycle.

%Now, we follow with the most difficult case.
%{\bf Case 4 ($\tau(G)> \alpha(G)+1 >2$).} 
By Lemma~\ref{conexo} we have that $q(H) \geq \alpha(H)-1+\Gamma(\alpha(H), \tau(H))$ 
whenever $H$ is a connected graph with $\tau(H)> \alpha(H)+1 >2$.
On the other hand, it is not difficult to see that if $H$ is a graph as in ${\rm (iv)}$, ${\rm (v)}$ or ${\rm (vi)}$, 
then $H$ is a $2$-connected graph with $q(H)=\alpha(H)+\Gamma(\alpha(H), \tau(H))$. 
Furthermore, since all that $q$-minimal graphs with $\tau(H)> \alpha(H)+1 >2$ are not $2$-connected, 
then $q(H) \geq \alpha(H)+\Gamma(\alpha(H), \tau(H))$
whenever $H$ is a $2$-connected graph with $\tau(H)> \alpha(H)+1 >2$.
%Hence,
%\[
%q(G) \leq \alpha(G)+\Gamma(\alpha(G), \tau(G)),
%\] 
%for all the $\{q,2\}$-minimal graphs $G$ with $\tau(G)\geq \alpha(G)+2$.

%\medskip

Therefore, to finish the proof we only need to 
show that $G$ is a $\{q,2\}$-minimal graph with $\tau(G)\geq \alpha(G)+2$ if and only if 
$G$ is as in   ${\rm (iv)}$, ${\rm (v)}$ or ${\rm (vi)}$. 
In order to do so, we  follow the same sequence of arguments 
as in the proof of Lemma~\ref{conexo}.
%Following the notation in Lemma~\ref{conexo}:

\medskip

Let $e\in E(G)$, if $e$ is not an $\alpha$-critical edge of $G$, then $G\setminus e$
is a $q$-minimal graph and $G$ is $2$-connexion of a PTT graph.
Therefore we can assume that $G$ is an $\alpha$-critical graph. 
%since
%if $G$ is not $\alpha$-critical, then $G\setminus e$ is a $q$-minimal graph
%for some edge $e$ and the result follows readily. 
\begin{Claim}\label{claim3} 
If %$G$ be a $q$-minimal and $\alpha$-critical graph and let 
$v$ is a vertex of $G$ of maximal degree, then either $G-v$ is $q$-minimal or $\{q,2\}$-minimal.
\end{Claim}
\demo
If $G\setminus v$ is neither  $q$-minimal nor $\{q,2\}$-minimal, 
then $q(G\setminus v)\geq \alpha(G)+\Gamma(\alpha(G), \tau(G)-1)+1$.
Therefore, the result follows using the same arguments that in Claim~\ref{claim1}.
\QED

%\begin{description}
%\item[{\it (1)}] $G-v$ is $q$-minimal or $\{q,2\}$-minimal,
%\item[{\it (2)}] 

%\item[{\it (4)}] 
%We already prove that, if $G$ is a $\{q,2\}$-minimal, 
%then $q(G) \leq \alpha(G)+\Gamma(\alpha(G), \tau(G))$.
%Furthermore, 
Since $G-v$ is either a $q$-minimal or a $\{q,2\}$-minimal graph, %(claim~\ref{claim3}), 
then by the induction hypothesis
\[
q(G-v) =
\begin{cases}
\alpha(G-v)-1+\Gamma(\alpha(G-v),\tau(G-v)) & \mbox{ if } G \mbox{ is } q \mbox{-minimal},\\
\alpha(G-v)+\Gamma(\alpha(G-v),\tau(G-v)) & \mbox{ if } G \mbox{ is } \{q,2\} \mbox{-minimal}.\\
\end{cases}
\]
Since $\deg(v)=q(G)-q(G-v)$ and $\Gamma(\alpha(G),\tau(G))-\Gamma(\alpha(G),\tau(G)-1)=\left \lfloor \frac{|V(G)|-1}{\alpha(G)}\right \rfloor$ (Lemma~\ref{lemaG} $(iii)$), then
\begin{equation} \label{eq7}
\deg(v) =%q(G)-q(G-v) = 
\begin{cases}
%1+\Gamma(\alpha(G),\tau(G))-\Gamma(\alpha(G),\tau(G)-1)=
\left\lfloor \frac{|V(G)|-1}{\alpha(G)}\right \rfloor+1 & \mbox{ if } G \mbox{ is } q \mbox{-minimal},\\
\\
%\Gamma(\alpha(G),\tau(G))-\Gamma(\alpha(G),\tau(G)-1)=
\left \lfloor \frac{|V(G)|-1}{\alpha(G)}\right \rfloor & \mbox{ if } G \mbox{ is } \{q,2\} \mbox{-minimal}.\\
\end{cases}
\end{equation}

\begin{Claim}\label{claim4} 
Let $H$ be a PTT graph with $H_1,\ldots,H_a$ connected components and let $L$ be a $2$-connexion of $H$. 
If $L'$ is a maximal induced subgraph of $L$ with $\alpha(L')=\alpha(L)-1$, then $L'$ is given as in $(i)$, $(ii)$, and $(iii)$ in Claim~\ref{claim2}.
\end{Claim}
\demo
Let $e$ be a connexion edge of $L$. Since $L$ is a $\{q,2\}$-minimal graph, then $L\setminus e$ is a $1$-connexion of $H$.
Hence, applying Claim~\ref{claim2} to $L\setminus e$ we get the result.
\QED
%\demo
%If $L'$ %is a maximal induced subgraph of a $2$-connexion of a PTT graph $L$ with $\alpha(L')=\alpha(L)-1$ and 
%is not given as in Claim~\ref{claim2}, then $V(e_i)\subset V(L')$ for all the connexion edges, $e_i$, of $L$.
%Thus, by the $2$-connectivity of $L$, $\alpha(L')=\alpha(L)$; a contradiction.  
%%Therefore, $L'$ is a maximal induced subgraph of a $2$-connexion of a PTT graph $L$ with $\alpha(L')=\alpha(L)-1$ 
%%if and only if $L'$ is given as in Claim~\ref{claim2}.  
%\QED
%Since $k=\lfloor \frac{|V(G-v)|}{\alpha(G)} \rfloor \leq |H_{i_0}|\leq \lceil \frac{|V(G-v)|}{\alpha(G)} \rceil=k+1$,
%then using the equation $(\ref{eq6})$ and Claim~\ref{claim2}, we obtain that
%${\rm deg}(v)\geq 2(k+1-2)+(k+1-3)$.
%%${\rm deg}(v)\geq 2(\lceil \frac{|V(G-v)|}{\alpha(G)} \rceil-2)+(\lceil \frac{|V(G-v)|}{\alpha(G)} \rceil-3)$.
%If $k+1 \geq 4$, then ${\rm deg}(v)\geq k+1+1=\lceil \frac{|V(G-v)|}{\alpha(G)} \rceil+1$;
%%If $\lceil \frac{|V(G-v)|}{\alpha(G)} \rceil \geq 4$, then
%%${\rm deg}(v)\geq \lceil \frac{|V(G-v)|}{\alpha(G)} \rceil+1$;
%a contradiction. %to equation~(\ref{eq7}).
%That is, not there exists $\{q,2\}$-minimal when $\lceil \frac{|V(G-v)|}{\alpha(G)} \rceil \geq 4$.
%Therefore we can assume that $\lceil \frac{|V(G-v)|}{\alpha(G)} \rceil\leq 3$.
%By claim~\ref{claim3} we only have two possibilities:
%
\medskip

Now, we will consider the cases when $G-v$ is either $q$-minimal or $\{q,2\}$-minimal:

\medskip

\noindent {\bf Case ($G-v$ is $q$-minimal).} 
Take $L=G\setminus v$ and $L'=G\setminus N[v]$.
If $G\setminus N[v]$ is as in Claim~\ref{claim2} ${\rm (i)}$, then $G$ must be a complete graph.
If $G\setminus N[v]$ is as in Claim~\ref{claim2} ${\rm (ii)}$, then $G$ must be an odd subdivision of the complete graph ${\cal K}_4$.

Now, assume that $G\setminus N[v]$ is as in Claim~\ref{claim2} ${\rm (iii)}$.
Using equations~(\ref{eq6}) and~(\ref{eq7}), we get that $k+1=\deg(v)\geq 2k-2$, where $k=\lfloor \frac{|V(G)|-1}{\alpha(G)} \rfloor$, that is, $k\leq 3$.
If $k=2$, then $G$ is either an odd cycle or an odd subdivision of ${\cal K}_4$ and if $k=3$, then $G$ is ${\cal H}$.
%\begin{itemize}
%\item $G-v$ has no odd cycle.
%Since $\lceil \frac{|V(G-v)|}{\alpha(G)} \rceil= 3$, 
%then $H_{i_0}$ is a complete graph with two or three vertices. 
%If $H_{i_0}={\cal K}_3$, then $|P|=2$ and $G_u$ 
%(the graph joined to $H_{i_0}$ by only one connexion edge 
%with one end equal to $u \in P$)
%is either an even path for all $u \in P$ 
%(since if the subjacent graph of $G_u$ has one component isomorphic 
%to ${\cal K}_3$ which is not a leaf of $G_u$, then $G$ will not be $\alpha$-critical)
%or equal to an even path with exactly one leaf isomorphic to ${\cal K}_3$. 
%The first possibility give us that $G$ is an even subdivision of 
%${\cal K}_4$ and in second possibility we have that by the $\{q,2\}$-minimality of $G$, it must be isomorphic to ${\cal H}$. 
%If $H_{i_0}={\cal K}_2$, then $G_u$ is an even path for all 
%$u \in P=V(H_{i-0})$ and therefore $G$ is an odd cycle.
%
%\item $G-v$ has an odd cycle.
%In this case, either the odd cycle has exactly three 
%vertices that are incident $v$ (by Claim~\ref{claim2}(ii))
%or the odd cycle has at least two vertices that are incident 
%to a connexion edge (by the $2$-connectivity of $G$). 
%Therefore using that $\deg(v)\leq 3$ (by equation~(\ref{eq7}))
%is not difficult to see that $G$ must be an even subdivision of ${\cal K}_4$.
%\end{itemize} 

\medskip
 
\noindent {\bf Case ($G-v$ is $\{q,2\}$-minimal).} 
Take $L=G\setminus v$ and $L'=G\setminus N[v]$.
If $L'$ is as in Claim~\ref{claim4} ${\rm (i)}$, then $G$ is a $2$-connexion of a PTT graph, but it is not an $\alpha$-critical graph.
If $G\setminus N[v]$ is as in Claim~\ref{claim2} ${\rm (ii)}$, then $G$ is an odd subdivision of ${\cal K}_4$.

Now, assume that $G\setminus N[v]$ is as in Claim~\ref{claim4} ${\rm (iii)}$.
Using equations~(\ref{eq6}) and~(\ref{eq7}), we get that $k=\deg(v)\geq 2k-2$, where $k=\lfloor \frac{|V(G)|-1}{\alpha(G)} \rfloor$, that is, $k\leq 2$.
If $k=2$, then $G$ is an odd cycle.

Finally, $G-v$ is not an odd subdivision of ${\cal K}_4$ because $\deg(v)\overset{\tiny (\ref{eq7})}{=}2$ and if
${\cal O}$ is an odd subdivision of ${\cal K}_4$, then $\alpha({\cal O}\setminus\{a,b\} )=\alpha({\cal O}) $ for all $a,b\in V({\cal O})$.
If $G-v$ is ${\cal H}$, then $\alpha({\cal H})=2$, $\omega({\cal H})=3$, and $\deg(v)\overset{\tiny (\ref{eq7})}{=}3\geq |V(G)|-\omega({\cal H})= 4$; a contradiction.
%By the induction hypothesis, there exists a vertex $u\neq v$ such that $G\setminus\{u,v\}$ is a $q$-minimal graph.
%Following the same arguments that we used to show that $G$ is a $q$-minimal graph we obtain the result.
%\vspace{-7.4mm}
\QED

%The methods introduced in this paper for the connected and $2$-connected 
%cases can be extended to the $3$-connected case and more generally to 
%the $k$-connected case. These will appear in subsequent papers.

\medskip 

\begin{Remark} 
{\rm After this paper was firstly submitted in 2006, the authors realized that 
Theorem~\ref{minimal} was also obtained independently in \cite{doignon}.} 
\end{Remark} 

\begin{Remark} 
{\rm It can be proved that for any fixed $\delta_{-}(G)=\alpha(G)-\tau(G)=k> 0$ 
there exist a finite number of ``basic" graphs such that if $G$ is $\{q,2\}$-minimal graph 
with $\delta_-(G)=k$, then $G$ is an odd subdivision of some of this basic graphs.
For instance, if $G$ is a $\{q,2\}$-minimal graph with $\delta_-(G)=1$, then $G$
is an odd subdivision of  the complete bipartite graph ${\cal K}_{2,3}$. 
%In this case, 
%The classification of the $\{q,2\}$-minimal graphs with $\tau(G)=\alpha(G)>1$ is more difficult. %to explain. 
%However, is easy to see that the $\{q,2\}$-minimal graphs with $\tau(G)=\alpha(G)>1$ are closed under odd subdivision.
} 
\end{Remark} 
%==================================================================% 
%==============================Bounds==============================% 
%==================================================================% 

\section{Some bounds for the stability and covering number of a graph}\label{bounds} 

The following results are in the spirit of \cite{bounds}, 
where the authors were motivated in bounding invariants for edge rings. 
In this paper, we concentrate mainly on the combinatorial aspects of these bounds. 

\medskip 

The theorem below gives an idea of the class of graphs that are Cohen-Macaulay and 
of those graphs that are far from being Cohen-Macaulay. 
We thank N. Alon (private communication) for some useful suggestions for making the 
proof of this result simpler and more readable.
 
\begin{Theorem}\label{alon} 
Let G be a graph without isolated vertices, then 
$$ 
\alpha(G)\leq \tau(G)[1+\delta(G)]. \vspace{-1mm}
$$ 
\end{Theorem} 
\demo 
First, let fix a minimal vertex cover $C$ with $\tau(G)$ vertices. %, where $k=\tau(G)$. 
Label the vertices of $C$ from $1$ to $\tau(G)$.
For each $i\in C$, let $T_i$ be a maximal stable set containing $i$,
with $|T_i|\geq \sigma_v(G)$.
Let $k$ be the minimal natural number such that
$$ 
C \subseteq \bigcup_{i=1}^{k} T_i.
$$
Clearly $0< k \leq \tau(G)$.
Let $M=V(G)\setminus C$ and take $C_i=C\cap T_i$ and $M_i=M\cap T_i$ for all $i=1,\ldots,\tau(G)$.
Since $M$ is a maximal stable set and $G$ does not have isolated
vertices, then for each vertex $v\in M$ there is an edge $e=\{v,v'\}$ with $v'\in C$.
That is, 
 
\begin{equation}\label{eq0} 
M = \bigcup_{i=1}^k (M \cap N(C_i)). 
\end{equation} 

\medskip 

Since $S_i=V(G)\setminus T_i=(C \setminus C_i) \cup (M \setminus M_i)$ 
is a minimal vertex cover with $|S_i|\leq n-\sigma_{v}(G)$ 
for all $i=1,\ldots,k$, then 
\[ 
|C \setminus C_i| + |M \setminus M_i|=|(C \setminus C_i) \cup (M \setminus M_i)|= |S_i| \leq n-\sigma_{v}(G). 
\] 
Hence, as $M \cap N(C_i) = M \setminus M_i$, then
\begin{equation}\label{eq1} 
\begin{array}{cl} 
|M \cap N(C_i)| = |M \setminus M_i|&  \leq n-\sigma_{v}(G)-|C \setminus C_i|\\ 
                                     &  = |C|+|M| -\sigma_{v}(G)-|C \setminus C_i|\\ 
                                     &  =|C_i|+\alpha(G)-\sigma_{v}(G) = |C_i|+\delta(G). 
\end{array} 
\end{equation} 

Taking 
\[ 
A_i=C_i\setminus (\bigcup_{j=1}^{i-1}C_j)
\mbox{ and } 
B_i=(M \cap N(C_i))\setminus (\bigcup_{j=1}^{i-1}M \cap N(C_j)), 
\] 
we have that 
\begin{equation}\label{eq2} 
|C_i\setminus A_i| \leq |M\cap N(C_i\setminus A_i)|.
\end{equation} 

Indeed, if $|C_i\setminus A_i| > |M\cap N(C_i\setminus A_i)|$, then 
$ 
C\setminus(C_i \setminus A_i) \cup (M\cap N(C_i\setminus A_i)) 
$ 
would be a vertex cover of cardinality $|C\setminus(C_i \setminus A_i)| + |M\cap N(C\setminus A_i)|<|C|$; a contradiction. 

\medskip 

To finish the proof, we use the inequalities (\ref{eq1}) and (\ref{eq2}) 
to conclude that 
\begin{equation}\label{eq3} 
\begin{array}{cl} 
|B_i|&=|M \cap N(C_i)|-|(M \cap N(C_i)) \cap ( \bigcup_{j=1}^{i-1}(M \cap N(C_j)))|\\ 
&=|M \cap N(C_i)|-|M \cap N(C_i) \cap  N(\bigcup_{j=1}^{i-1}C_j))|\\ 
&\overset{(\ref{eq1})}{\leq} |C_i|+\delta(G)-|M \cap N(C_i \cap \bigcup_{j=1}^{i-1}C_j))|\\ 
&\overset{(\ref{eq2})}{\leq}  |C_i|+\delta(G)-|C_i\setminus A_i|=|A_i|+\delta(G). 
\end{array} 
\end{equation}

Therefore 
\begin{equation}\label{eq5} 
\begin{array}{cl} 
\alpha(G)=|M| & \overset{(\ref{eq0})}{=}|\bigcup_{i=1}^k (M \cap N(C_i))|
=\sum_{i=1}^k |B_i| \overset{(\ref{eq3})}{\leq} \sum_{i=1}^k (|A_i|+ \delta(G))\\ 
\\ 
&=\sum_{i=1}^k |A_i|+ \sum_{i=1}^k \delta(G) \leq |C|+ \tau(G)\delta(G)
= \tau(G)[1+\delta(G) ] \vspace{-7mm} 
\end{array} 
\end{equation} \QED 

\medskip

When $\delta(G)>0$ it is not difficult to characterize the graphs $G$ such that $\alpha(G)= \tau(G)[1+\delta(G)]$.
Let $\tau$ and $\delta$ be positive numbers and let $H^-(\tau,\delta)$ and $H^+(\tau,\delta)$ be the graphs 

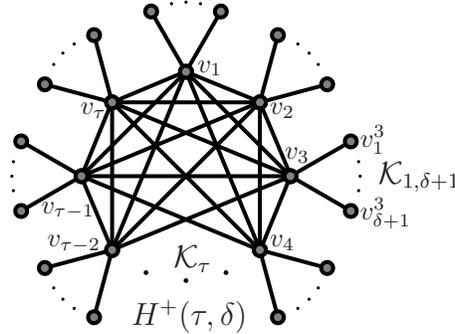
\vspace{35mm} 
\begin{figure}[h] 
\begin{center} 
$ 
\setlength{\unitlength}{.4mm} 
\thicklines\begin{picture}(125,0) 

\scalebox{.75}{\includegraphics{figura4e}} 

{\large 
\put(-63,20){${\cal K}_{\tau}$} 
\put(-77,1){$H^+(\tau,\delta)$}
\put(5,48){${\cal K}_{1,\delta+1}$} 
}
{\small
\put(-55,85){$v_1$}
\put(-31,71){$v_2$}
\put(-94,71){$v_{\tau}$} 
\put(-26,56){$v_3$}
\put(-107,38){$v_{\tau-1}$}
\put(-31,27){$v_4$}
\put(-105,27){$v_{\tau-2}$} 
\put(-1,60){$v^3_1$}
\put(-1,36){$v^3_{\delta+1}$}
}  
\end{picture} 
$ 
\end{center} 
\end{figure} 
%\vspace{-9mm}
%\mbox{.}

\noindent on the vertex set $V(H^-(\tau,\delta))=V(H^+(\tau,\delta))=V_0\cup V_1 \cup \cdots \cup V_{\tau} $ 
where $V_0=\{v_1,v_2,\ldots,v_{\tau}\}$, $V_i=\{v^i_1,\ldots, v^i_{\delta+1}\}$ for all $i=1,\ldots,\tau$ 
and edge sets 
\[
E(H^-(\tau,\delta))=(\bigcup_{i=1}^{\tau} \{ \{v_i,v^i_j \} \, | \, 1\leq j \leq\delta+1\})
\] 
and 
\[
E(H^+(\tau,\delta))= E(H^-(\tau,\delta))\cup \{\{v_i,v_j\} \, | \, 1 \leq  i \neq j \leq \tau\}.
\]

\medskip

A set of edges in a graph $G$ is called independent or a {\it matching} if no two of them have a vertex in common. 
A pairing by an independent set of edges of all the vertices of a graph $G$ is called a {\it perfect matching}. 

\begin{Corollary}
Let G be a graph without isolated vertices. 
\begin{description}
\item[{\it (i)}] If $\delta(G)> 0$, then $\alpha(G)= \tau(G)[1+\delta(G)]$ if and only if 
\[
E(H^-(\tau,\delta)) \subseteq E(G) \subseteq E(H^+(\tau,\delta)),
\] 
where $\tau=\tau(G)$ and $\delta=\delta(G)$.

\item[{\it (ii)}] If $\delta(G)=0$ and $\alpha(G)= \tau(G)$, then $G$ has a perfect matching. 
\end{description}
\end{Corollary}
\demo
We  use the same notation as in the proof of Theorem~\ref{alon}.

{\bf (i)}
%Since $\delta(G)>0$ and $\alpha(G)= \tau(G)[1+\delta(G)]$, 
%then 
Since $\delta(G)>0$ and $\alpha(G)= \tau(G)[1+\delta(G)]$, then using equation~(\ref{eq5}) %in the proof of Theorem~\ref{alon} 
we can conclude that $k=\tau(G)$.

Following the proof of Theorem~\ref{alon}, we have that $|C\cap M'|\leq1$ for all 
$M'$ maximal stable set. %with $|M'|\geq \sigma_v(G)$.
Moreover, for all $u\in C$ there exists a $M'$ maximal stable set  with %$|M'|\geq \sigma_v(G)$ and 
$C\cap M'=\{u\}$.
Thus, the equation~(\ref{eq1}) reduces to $|M\cap N(u)|\leq 1+\delta(G)$ for all $u\in C$, where $M=V(G)\setminus C$.
On the other hand, since $M=\cup_{u\in C}(M\cap N(u))$ and $\alpha(G)= \tau(G)[1+\delta(G)]$, 
then $|M\cap N(u)|=1+\delta(G)$ for all $u\in C$ and
 $(M\cap N(v))\cap (M\cap N(u))=\emptyset$ for all $u\neq v\in C$.
Furthermore, since $M$ is a stable set, then $E(H^-(\tau(G),\delta(G)))\subseteq E(G)\subseteq E(H^+(\tau(G),\delta(G)))$.

Finally, note that if $E(H^-(\tau,\delta))\subseteq E(G)\subseteq E(H^+(\tau,\delta))$ for some $\tau>0$ and $\delta>0$, 
then clearly $\alpha(G)= \tau(G)[1+\delta(G)]$.

{\bf (ii)} Following the proof of Theorem~\ref{alon} we have that
$(ii)$ reduces to prove that for all $i=1,\ldots,k$ the induced subgraph $G_i=G[A_i\cup B_i]$ has a perfect matching, 
that is, $\nu(G_i)=|A_i|=|B_i|$ for all $i=1,\ldots,k$. 

Since $G_i$ is a bipartite graph ($A_i$ and $B_i$ are stable sets of $G$), then by Konig's theorem $\nu(G_i)=\tau(G_i)$.
Hence, we only need to prove that $|A_i|=|B_i|$ and $\tau(G_i)=|A_i|$ for all $i=1,\ldots,k$.

Firstly, since  $C=\sqcup^k_{i=1} A_i$ and $M=\sqcup^k_{i=1} B_i$, then $\sum_{i=1}^k |A_i|=\tau(G)=\alpha(G)=\sum_{i=1}^k |B_i|$.
On the other hand, since $\delta(G)=0$, then the equation~(\ref{eq3}) reduces to $|A_i|\leq |B_i|$ and therefore $|A_i|=|B_i|$ for all $i=1,\ldots,k$.

Finally, we will prove that $\tau(G_i)=|A_i|$ for all $i=1,\ldots,k$.
Since $A_i$ is a vertex cover of $G_i$, then $\tau(G_i)\leq|A_i|$.
Furthermore, if $\tau(G_i)<|A_i|$, then there exist a stable set  $N$ of $G_i$ with $|N|>|A_i|$.
Since $M \cap N(\cup_{j=1}^i C_j)=\cup_{j=1}^i (M\cap N(C_i))=  \cup_{j=1}^i M \setminus M_j$, 
then $N\cup (T_i \setminus A_i)=(N\cap A_i)\cup M_i \cup (C_i\setminus A_i)\cup (N\cap B_i)\subset T_i\cup (N\cap B_i)$ 
is a stable set of ($N(C_i\setminus A_i)\cap B_i=\emptyset$) of $G$ with $|T_i|-|A_i|+|N|> |T_i|=\alpha(G)$ vertices; a contradiction.
\mbox{}\QED

%==================================================================% 
%==============================B-graphs=============================% 
%==================================================================% 

\subsection{$B$-graphs} \label{bgraph} 

A graph is called a {\it $B$-graph} if every vertex belongs to a maximum stable set (that is, to a stable set of largest size). 
This concept was introduced by Berge in~\cite{berge}. 

The {\it $\sigma_{v}$-cover} number of a graph, denoted by $\sigma_{v}(G)$, 
is the maximum natural number $m$, such that every vertex of 
$G$ belongs to a maximal independent set with at least $m$ vertices. 
%We define $\delta(G) = \alpha(G)-\sigma_{v}(G)$. %and $\omega_{v}(G)$ as $\sigma_{v}(\overline{G})$. 
Clearly, $G$ is a $B$-graph if and only if $\alpha(G)=\sigma_v(G)$ if and only if $\delta(G)=0$, where $\delta(G) = \alpha(G)-\sigma_{v}(G)$.
%We say that a graph is {\it $\tau$-critical} if $\tau(G-v) < \tau(G)$ for all the vertices $v\in V(G)$. 

\medskip

Now, we  define two %new important 
invariants that measure when a graph is a $B$-graph or a $\tau$-critical graph. 
Let
\[ 
\alpha_{core}(G)\, \, =\bigcap_{|{\cal M}_i|=\alpha(G)}^{stable \, set}{\cal M}_i \, \, \mbox{  and  }\, \, 
\tau_{core}(G)\, \, =\bigcap_{|{\cal C}_i|=\tau(G)}^{vertex \, cover}{\cal C}_i, 
\] 
be the intersection of all the maximum stable sets and of all the minimum vertex covers of $G$, respectively. 
Also, let $B_{\alpha\cap \tau}=V(G)\setminus (\alpha_{core}(G) \cup \tau_{core}(G))$. 

\begin{Example}\label{exa1}\rm 
To illustrate the concepts of $\alpha_{core}(G)$, $\tau_{core}(G)$ and  $B_{\alpha\cap \tau}$ consider the following graph: 

\vspace{24mm} 
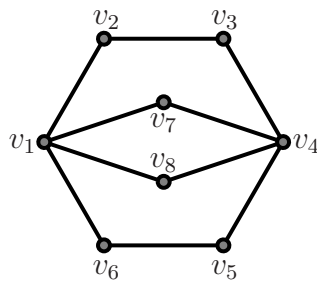
\begin{center} 
$ 
\setlength{\unitlength}{.28mm} 
\thicklines\begin{picture}(127,0) 

\scalebox{0.75}{\includegraphics{figura3e}} 

{\large 
\put(-135,51){$v_1$} 
\put(-95,110){$v_2$}
\put(-39,110){$v_3$}
\put(0,51){$v_4$} 
\put(-39,-9){$v_5$} 
 \put(-95,-9){$v_6$} 
\put(-68,61){$v_7$} 
\put(-68,42){$v_8$} 
} 
\end{picture} 
$ 
\end{center} 
\vspace{1mm} 
since $\alpha(G)=4$, $\tau(G)=4$ and 
$ 
\{v_2,v_5,v_7,v_8\}, \{v_2,v_6,v_7,v_8\}, \{v_3,v_5,v_7,v_8\}, \{v_3,v_6,v_7,v_8\}
$ 
are the maximum stable sets of G, then
\begin{itemize} 
\item $\alpha_{core}(G)=\{v_7,v_8\}$, 

\item $\tau_{core}(G)=\{v_1,v_4\}$, and

\item $B_{\alpha\cap \tau}=\{v_2,v_3,v_5,v_6\}$. 
\end{itemize} 

\end{Example} 

\medskip 

Since $M$ is a maximum stable set of $G$ if and only if $C=V(G)\setminus M$ is a minimum vertex cover of $G$, then
$G$ is a $B$-graph if and only if %every vertex of $G$ belongs to a maximum stable set if and only if 
$V(G)=\bigcup_{|{\cal M}_i|=\alpha(G)}^{stable \, set}{\cal M}_i $ if and only if $\tau_{core}(G)=\bigcap_{|{\cal C}_i|=\tau(G)}^{vertex \, cover}{\cal C}_i=\emptyset$. 
Similarly,  since a graph is {\it $\tau$-critical} if and only if $\tau(G-v) < \tau(G)$ for all $v\in V(G)$ if and only if there exists a maximum stable set $M_v$ of $G$ such that $v\notin M_v$ for all $v\in v(G)$,
then $G$ is a $\tau$-critical graph if and only if $\alpha_{core}(G)=\emptyset$. 

\begin{Proposition}\label{decomposition} 
Let $G$ be a graph, then 
\[ 
V(G)=\alpha_{core}(G) \sqcup \tau_{core}(G) \sqcup B_{\alpha\cap \tau}, 
\] 
furthermore 
\begin{description} 
\item[{\it (i)}] $G[\alpha_{core}(G)]$ is a trivial graph, 

\item[{\it (ii)}] $N(\alpha_{core}(G))\subseteq\tau_{core}(G)$, 

\item[{\it (iii)}] $G[B_{\alpha\cap \tau}]$ is both a $\tau$-critical graph 
as well as a $B$-graph without isolated vertices, and% and without isolated vertices graph. 

\item[{\it (iv)}] $\alpha(G)-|\alpha_{core}(G)| \leq \tau(G)-|\tau_{core}(G)|$
\end{description} 
\end{Proposition} 
\demo Firstly, is clear that $\alpha_{core}(G) \cap \tau_{core}(G)=\emptyset$.
Also, by the definition of $B_{\alpha\cap \tau}$ it is clear that 
$\alpha_{core}(G) \cap B_{\alpha\cap \tau}=\emptyset$ and $\tau_{core}(G) \cap B_{\alpha\cap \tau}=\emptyset$.

{\bf (i)}  Since $\alpha_{core}(G)$ is the intersection of stable sets, then $\alpha_{core}(G)$ is a stable set and therefore $G[\alpha_{core}(G)]$ is a trivial graph.

{\bf (ii)} Since %$G[V(G)\setminus \tau_{core}(G)]$ is a $B$-graph, then
$\tau_{core}(G) = V(G)\setminus \bigcup_{|{\cal M}_i|=\alpha(G)}^{stable \, set}{\cal M}_i $ and $\alpha_{core}(G) \subset V(G)\setminus \tau_{core}(G)$, then 
$\alpha_{core}(G)$ is the set of isolated vertices of $G[V(G)\setminus \tau_{core}(G)]$.
Therefore $N(\alpha_{core}(G))\subseteq\tau_{core}(G)$.

{\bf (iii)} Since $\alpha(G[B_{\alpha\cap \tau}])=\alpha(G)-|\alpha_{core}(G)|$, $\tau(G[B_{\alpha\cap \tau}])=\tau(G)-|\tau_{core}(G)|$ 
and $B_{\alpha\cap \tau}=V(G)\setminus (\alpha_{core}(G) \cup \tau_{core}(G))$, then $\alpha_{core} (B_{\alpha\cap \tau})=\emptyset$ 
and  $\tau_{core} (B_{\alpha\cap \tau})=\emptyset$.
%Using the observations previous this proposition, then 
Therefore $G[B_{\alpha\cap \tau}]$ is a $\tau$-critical graph and a $B$-graph without isolated vertices.

{\bf (iv)} Since $G[B_{\alpha\cap \tau}]$ is a $B$-graph, then $\delta(G[B_{\alpha\cap \tau}])=0$.
Therefore, applying Theorem~\ref{alon} to $G[B_{\alpha\cap \tau}]$, %we obtain that 
\[
\alpha(G)-|\alpha_{core}(G)|= \alpha(G[B_{\alpha\cap \tau}]) \leq \tau(G[B_{\alpha\cap \tau}])=\tau(G)-|\tau_{core}(G)|. \vspace{-8mm} 
\] 
\QED 

%Note that in example~\ref{exa1} we have that $\alpha(G)-|\alpha_{core}(G)|=4-2 = \tau(G)-|\tau_{core}(G)|$.

\begin{Remark}\rm 
If $v$ is an isolated vertex, then $v\in \alpha_{core}(G)$, and if ${\rm deg}(v) > \tau(G)$, 
then $v$ does not belong to any stable set with $\alpha(G)$ vertices and therefore $v\in \tau_{core}(G)$. 
Note that in general the induced graph $G[B_{\alpha\cap \tau}]$ is not necessarily connected.
\end{Remark} 

\begin{Corollary}{\rm (\cite[Proposition 7]{berge})}\label{coro2} 
If $G$ is a $B$-graph without isolated vertices, then $G$ is a $\tau$-critical graph. 
\end{Corollary} 
\demo 
Since $G$ is a $B$-graph, then $\tau_{core}(G)=\emptyset$. 
Thus, by Proposition~\ref{decomposition} (ii), $N(\alpha_{core}(G))=\tau_{core}(G)=\emptyset$. 
Moreover, since $G$ has no isolated vertices, then $\alpha_{core}(G)=\emptyset$. 
Therefore $G$ is a $\tau$-critical graph.   
\QED 

\begin{Remark}\rm 
The bound of Proposition~\ref{decomposition} {\rm (iv)}  improves 
the bound given in \cite[Theorem 2.11]{levit} for the number of vertices in $\alpha_{core}(G)$. 

Their result states that if $G$ is a graph of order $n$ and 
\[\alpha(G) > \frac{n+k-{\rm min} \{1,|N(\alpha_{core}(G))|\}}{2}, \mbox{ for some } k\geq 1, \]
then $|\alpha_{core}(G)| \geq k+1$. 
Moreover, if $n+k-{\rm min} \{1,|N(\alpha_{core}(G))|\}$ is even, then $|\alpha_{core}(G)|\geq k+2$. 

\medskip 

%On the other hand, 
Our result %(Proposition~\ref{decomposition} $(iv)$) 
states that if $G$ is a graph of order $n$ and  $\alpha(G)\geq \frac{n+k'}{2}$ for some $k'\geq 0$, then \vspace{-3mm}
\[
%\begin{array}{cl}
|\alpha_{core}(G)|  \overset{\rm (iv)}{\geq} \alpha(G)-\tau(G)+|\tau_{core}(G)| = 2\alpha(G)-n +|\tau_{core}(G)| \geq k' + |\tau_{core}(G)|.
%\end{array}
\]

In order to compare both bounds we can write their bound in the following equivalent way:
If $G$ is a graph of order $n$, $|N(\alpha_{core}(G))|=0$ ($|N(\alpha_{core}(G))|\geq 1$) and 
\[
\alpha(G) \geq
\begin{cases}
\frac{n+(k+1)}{2} \, (\frac{n+(k+1)}{2}) \mbox{ if } n+k \mbox{ is odd},\\
\\
\frac{n+(k+2)}{2} \, (\frac{n+k}{2}) \mbox{ if } n+k \mbox{ is even}.
\end{cases}
\]
\vspace{-5mm}
for some $k\geq 1$, then 
\[
|\alpha_{core}(G)| \geq
\begin{cases}
k+1 \, (k+2) \mbox{ if } n+k  \mbox{ is odd},\\
k+2 \, (k+1) \mbox{ if } n+k  \mbox{ is even}.
\end{cases}
\]

Since $|N(\alpha_{core}(G))|  \leq |\tau_{core}(G)|$ (Proposition~\ref{decomposition} {\rm (ii)}), then our bound improves their bound.
Furthermore, the bounds are equivalent if and only if $|N(\alpha_{core}(G))|  = |\tau_{core}(G)|=0,1$. %$k=k'\geq 2$
\end{Remark} 

%==================================================================% 
%============================Conjectures===========================% 
%==================================================================% 

\subsubsection{Conjectures} 

In this section we present a conjecture that generalizes the result obtained from 
Theorem~\ref{alon} when $G$ is a $B$-graph.
Before stating the conjecture we will introduce a new graph invariant. %that is a strong version of the 
%$\sigma_v$-cover number of the complement of a graph. %introduced in the beginning of this subsection. 

\begin{Definition}\rm 
The {\it $\omega_{e}$-clique covering} number of $G$, denoted by 
$\omega_{e}(G)$, is the greatest natural number $m$ so that every edge in $G$ 
belongs to a clique of size at least $m$. 
\end{Definition} 

\begin{Example}\label{exa2}\rm
In order to illustrate the previous concept consider the following graph: 

\vspace{34mm} 
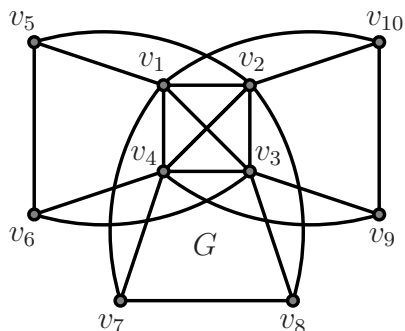
\begin{center} 
$ 
\setlength{\unitlength}{.36mm} 
\thicklines\begin{picture}(140,0) 

\scalebox{.65}{\includegraphics{figura1e}} 

{\large 
%\put(-72,100){$G$} 
\put(-72,20){$G$} 
\put(-92,90){$v_1$} 
\put(-55,90){$v_2$} 
\put(-49,55){$v_3$} 
\put(-95,55){$v_4$} 
\put(-140,105){$v_5$} 
\put(-140,25){$v_6$} 
\put(-107,-7){$v_7$} 
\put(-40,-7){$v_8$} 
\put(-8,105){$v_{10}$} 
\put(-7,25){$v_9$} 
} 
\end{picture} 
$ 
\end{center} 
\vspace{-2mm} 

For this graph we have that: 
\begin{itemize} 
\item $\omega(G)=\alpha(\overline{G})=4$ because $\{v_1,v_2,v_3,v_4\}$ is a clique of $G$,

\item $\omega_{e}(G)=2$ because the edge $\{v_5,v_6\}$ is not in any induced ${\cal K}_3$ of $G$, 

\item %$\sigma_{e}(G)=
$\alpha(G)=\omega_{e}(\overline{G})=4$ because $\{v_1,v_6,v_8,v_9\}$, $\{v_2,v_6,v_7,v_9\}$ 
$\{v_3,v_5,v_7,v_{10}\}$ and $\{v_4,v_5,v_8,v_{10}\}$ are stable sets of $G$, and

\item %$\omega_{v}(G)=
$\sigma_{v}(\overline{G})=3$ because $\{v_{i-1},v_{i},v_{10-2i}\}$ for $i=1,2,3$ 
and $\{v_{1},v_{2},v_{5}\}$, $\{v_{1},v_{4},v_{7}\}$, $\{v_{3},v_{4},v_{9}\}$ are cliques of $G$. 

\end{itemize} 
\end{Example}

%\medskip
%Also, we define $\sigma_{e}(G)$ as $\omega_{e}(\overline{G})$ and $\omega_{v}(G)$ as $\sigma_{v}(\overline{G})$. 
%Is not difficult to see that %this invariants %of a graph $G$ (or hypergraph) satisfy that
The {\it $\omega_{e}$-clique covering} and the {\it $\sigma_{v}$-cover} numbers of $G$ satisfy the following two identities:
$$ 
\omega_{e}(G) \leq \sigma_{v}(\overline{G}) \leq \omega(G) \mbox{ and } \omega_{e}(\overline{G}) \leq \sigma_{v}(G) \leq \alpha(G). 
$$ 

\begin{Conjecture}\label{conje-prin} 
Let $G$ be a $B$-graph ($\sigma_{v}(G)=\alpha(G)$) without isolated vertices, then 
\[ 
\omega_{e}(G)\sigma_v(G)\leq |V(G)|. 
\] 

Furthermore, for all maximum stable sets $M$, there exist disjoint sets 
\[ 
A_{j} \subset V(G) \mbox{ for all } j=1,\ldots,|M|, 
\] 
such that 
\begin{description} 
\item[{\it (i)}] $|M\cap A_{j}|=1$ for all $j=1,\ldots,|M|$, 
\item[{\it (ii)}] $G[A_j]$ is a clique of order $\omega_{e}(G)$ for all $i=1,\ldots,|M|$. 
\end{description} 
\end{Conjecture} 

%%%%%%In order to illustrate this conjecture is not difficult to check that graph in the example~\ref{exa2} satisfy the conjecture.

%Note that, we cannot weaken the Conjecture~\ref{conje-prin}, for example, 
%we cannot change the invariant $\omega_{e}(G)$ by the invariant $\omega_{v}(G)$ in the formula of~\ref{conje-prin}. 
%To see this, consider the following graph: 

%However 
%$ 
%\omega_{v}(G)\alpha(G)=(3)(4)>10=n. 
%$

\begin{Remark}\rm
If $G$ is a graph without isolated vertices, then $\omega_{e}(G)\geq 2$.
Thus, if Conjecture~\ref{conje-prin} holds, then $2\alpha(G)\leq |V(G)|$  when $G$ is a $B$-graph without isolated vertices.
On the other hand, Theorem~\ref{alon} implies that if $G$ is a $B$-graph without isolated vertices, then $\alpha(G)\leq \tau(G)$.
%($\tau(G)=|V(G)|-\alpha(G)\geq \frac{|V(G)|}{2}\geq \alpha(G)$)
Since $\tau(G)=|V(G)|-\alpha(G)$, then
Conjecture~\ref{conje-prin} implies the bound given in Theorem~\ref{alon} when $G$ is a $B$-graph.
\end{Remark}

\begin{Remark}\rm
A weaker version of Conjecture~\ref{conje-prin} is given in \cite{galvin} and \cite{goddard}. 
In these papers the authors prove that, if $G$ is a graph with $n$ vertices such that every vertex 
belongs to a clique of cardinality $q+1$ and a stable set of cardinality $p+1$, then $|V(G)| \geq p+q+\sqrt{4pq}$.  
Using our terminology, this bound can be writen as 
\[
\sigma_v(G)+\sigma_v(\overline{G})-2+2\sqrt{(\sigma_v(G)-1)(\sigma_v(\overline{G})-1)} \leq |V(G)|.
\]
Since
$
\sigma_v(G)+\sigma_v(\overline{G})-2+2\sqrt{(\sigma_v(G)-1)(\sigma_v(\overline{G})-1)} \leq {\rm max}\{4(\sigma_v(G)-1), 4(\sigma_v(\overline{G})-1) \}
$, 
then in comparison with our bound, this bound is a bad lower bound for the number of vertices of a $B$-graph $G$.
%Also, note that this bound show that if $G$ is a $B$-graph with $\omega_e(G)\leq 3$, then our bound is true.
\end{Remark}

\begin{Remark}\rm 
If we do not assume that $G$ is a $B$-graph, then Conjecture~\ref{conje-prin} is false.
To see this fact, consider the following graph: 

\vspace{35mm} 
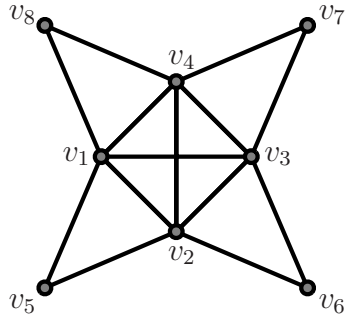
\begin{center} 
$ 
\setlength{\unitlength}{.40mm} 
\thicklines\begin{picture}(89,0) 

\scalebox{.80}{\includegraphics{figura2e}} 

{\large 
\put(-85,46){$v_1$} 
\put(-50,13){$v_2$} 
\put(-18,46){$v_3$} 
\put(-50,78){$v_4$} 
\put(-103,-3){$v_5$} 
\put(0,-3){$v_6$} 
\put(0,93){$v_7$} 
\put(-103,93){$v_8$} 
} 
\end{picture} 
$ 
\end{center} 
\vspace{0mm} 

For this graph we have, 
\begin{itemize} 
\item $\alpha(G)=4$ because $\{v_5,v_6,v_7,v_8\}$ is a stable set, 
\item $\sigma_v(G)=3$ because $\{v_1,v_6,v_7\}$, $\{v_2,v_7,v_8\}$, $\{v_3,v_5,v_8\}$, and $\{v_4,v_5,v_6\}$ are stable sets,
%\item $\omega(G)=4$ because $\{v_1,v_2,v_3,v_4\}$ is a clique,
\item $\omega_{e}(G)=3$ because all the edges are in an induced ${\cal K}_3$,
\item  $\omega(G)=4$, $\sigma_{e}(G)=3$ and $\omega_v(G)=3$ because $\overline{G}\cong G$, and
\item $G$ is not a $B$-graph because $\sigma_v(G)=3 \neq 4=\alpha(G)$.
\end{itemize} 

However, 
$$ 
\omega_{e}(G)\sigma_{v}(G)=(3)(3)=9> 8=n. 
$$ 
\end{Remark}

\paragraph{Hypergraphs}

A {\it hypergraph} ${\cal H}$ is a pair ${\cal H} =(V,{\cal E})$ where $V$ is a set of elements, called vertices, and ${\cal E}$ is a set of non-empty subsets of $V$ called hyperedges.
We say that a hypergraph ${\cal H}=(V,{\cal E})$ is a {\it $r$-uniform} hypergraph if  $|E|=r$ for all $E\in {\cal E}$.
A vertex $v\in V$ of a hypergraph ${\cal H}=(V,{\cal E})$ is called isolated if $v\notin E$ for all $E\in {\cal E}$.

\medskip

A subset $M$ of vertices of ${\cal H}$ is called a {\it stable} set if no two vertices in $M$ belong to a hyperedge of ${\cal H}$.  
We say that $M$ is a {\it maximal stable} set if it is maximal with respect to inclusion. 
The {\it stability number} of a hypergraph ${\cal H}$ is given by 
$$ 
\alpha({\cal H})={\rm max}\{|M| \, | \, M\subset V({\cal H}) \mbox{ is a stable set in } {\cal H}\}. 
$$ 

The {\it $\sigma_{v}$-cover} number of a hypergraph ${\cal H}$, denoted by $\sigma_{v}({\cal H})$, 
is the maximum natural number $m$ such that every vertex of 
${\cal H}$ belongs to a maximal independent set of ${\cal H}$ with at least $m$ vertices.

%Let ${\cal H}=(V,{\cal E})$ be a $r$-uniform hypergraph,
%a subset $W$ of $V$ is called a {\it clique} if any $r$ vertices in $W$ are a hyperegde of ${\cal H}$. 
%We call $W$ {\it maximal} if it is maximal with respect to inclusion. 
%The {\it clique number} of a hypergraph ${\cal H}$ is given by 
%$$ 
%\omega({\cal H})={\rm max}\{ |W| \, | \, W\subset V({\cal H}) \mbox{ is a clique in }  {\cal H}\}. 
%$$

%In a similar way that for graphs, we can define the $\sigma_v$-cover number, denoted by $\sigma_v({\cal H})$, 
%and the $\omega_e$-clique covering number, denoted by $\omega_e({\cal H})$, of a hypergraph ${\cal H}$.

\medskip

The next conjecture was stated in~\cite[ Conjecture 3.2.12]{thesis}. 
\begin{Conjecture}\label{consecuencia2} 
Let ${\cal H}=(V,{\cal E})$ be a $r$-uniform hypergraph without isolated vertices. %with $|V|=n$. 
%and with all the hyperedges of size $r$. 
If $\sigma_{v}({\cal H})=\alpha({\cal H})$, then 
$$ 
r\sigma_{v}({\cal H}) \leq |V|. 
$$ 
\end{Conjecture} 

The last conjecture follows from Conjecture~\ref{conje-prin} by the following argument:
Let ${\cal H}$ be a hypergraph and consider the graph $G({\cal H})$ 
defined on the same vertex set of ${\cal H}$ and for which $v_1,v_2 \in G({\cal H})$ are 
adjacent if and only if they are adjacent in ${\cal H}$. 

Clearly $G({\cal H})$ has the same stability number as ${\cal H}$.
Also, observe that $G({\cal H})$ and ${\cal H}$ have the same $\sigma_{v}$-cover number.  
Moreover, if ${\cal H}$ is a $r$-uniform hypergraph, then $r \leq \omega_e(G({\cal H}))$. %is at least or equal to the $\omega_e$-clique covering
%number of $G({\cal H})$. 
 %clique number of $G({\cal H})$ is greatest or equal to the greatest cardinality of an edge in ${\cal H}$. 
%Moreover, $G({\cal H})$ and ${\cal H}$ have the same $\sigma_{v}$ and $\omega_{e}$ covering numbers. 
Applying Conjecture~\ref{conje-prin} to the graph $G({\cal H})$, we have that, 
\[ 
r\sigma_{v}({\cal H}) \leq \omega_{e}(G({\cal H}))\sigma_{v}(G({\cal H})) \leq |V|.
\] 
%\medskip 
%The following conjecture differs from the ones 
%given above, in that it stresses the symmetry of the formulas 
%with respect to the complement of a graph. 

%\begin{Conjecture}\label{simetrica1} 
%Let $G$ be a graph with $n$ vertices and without isolated vertices. 

%\medskip 

%If $\alpha(G)=\sigma_{e}(G)$ and $\omega(G)=\omega_{e}(G)$, then 
%\[ 
%\omega_{e}(G)\sigma_{e}(G)\leq n. 
%\] 
%\end{Conjecture}

\noindent {\bf Acknowledgment}

The authors gratefully acknowledge the valuable and helpful comments pointed out by the
anonymous referee.

%The authors thank the referee, who contributed to the paper by pointing out some errors
%and making it more readable.
%\end{Ackno}
%==================================================================% 
%===========================Bibliografia===========================% 
%==================================================================% 

\end{document}